\documentclass[11pt]{article}
\usepackage[authoryear]{natbib}
\usepackage{latexsym, epsfig, amssymb, amsmath, amsthm, graphicx, mathrsfs}
\usepackage{anysize}
\usepackage[OT1]{fontenc}
\usepackage{hypernat}

\renewcommand{\baselinestretch}{1.1}
\marginsize{1.2in}{1in}{0.55in}{1.55in} %
\makeatletter
\def\singlespace{\def\baselinestretch{1}\@normalsize}

\makeatother

\newcommand{\bff}{\mbox{\bf f}}

\newcommand{\bu}{\mbox{\bf u}}
\newcommand{\bv}{\mbox{\bf v}}

\newcommand{\bx}{\mbox{\bf x}}
\newcommand{\by}{\mbox{\bf y}}

\newcommand{\bA}{\mbox{\bf A}}

\newcommand{\bD}{\mbox{\bf D}}

\newcommand{\bG}{\mbox{\bf G}}

\newcommand{\bI}{\mbox{\bf I}}

\newcommand{\bR}{\mbox{\bf R}}

\newcommand{\bX}{\mbox{\bf X}}

\newcommand{\bZ}{\mbox{\bf Z}}

\newcommand{\bzero}{\mbox{\bf 0}}
\newcommand{\bveps}{\mbox{\boldmath $\varepsilon$}}
\newcommand{\balpha}{\mbox{\boldmath $\alpha$}}
\newcommand{\bbeta}{\mbox{\boldmath $\beta$}}

\newcommand{\btheta}{\mbox{\boldmath $\theta$}}

\newcommand{\bmu}{\mbox{\boldmath $\mu$}}

\newcommand{\bomega}{\mbox{\boldmath $\omega$}}
\newcommand{\bSig}{\mbox{\boldmath $\Sigma$}}

\newcommand{\hbbeta}{\widehat\bbeta}
\newcommand{\hdelta}{\hat{\delta}}

\newcommand{\hbtheta}{\widehat\btheta}
\newcommand{\hbeta}{\widehat\beta}

\newcommand{\calA}{{\cal A}}
\newcommand{\sbbeta}{\mbox{\scriptsize \boldmath $\beta$}}
\newcommand{\sM}{{\scriptsize \cal M}}
\newcommand{\sA}{{\scriptsize \cal A}}

\newcommand{\var}{\mathrm{var}}
\newcommand{\cov}{\mathrm{cov}}
\newcommand{\corr}{\mathrm{corr}}
\newcommand{\Sig}{\mathbf{\Sigma}}

\newcommand{\diag}{\mathrm{diag}}

\newcommand{\sgn}{\mathrm{sgn}}
\newcommand{\supp}{\mathrm{supp}}

\newcommand{\fair}{\hdelta_{\text{FAIR}}}

\def\t{^T}
\def\toD{\overset{\mathscr{D}}{\longrightarrow}}
\def\toP{\overset{\mathrm{P}}{\longrightarrow}}

\begin{document}

\title{Invited Review Article: A Selective Overview of Variable Selection in High Dimensional Feature Space %
\date{September 1, 2009}
\author{Jianqing Fan and Jinchi Lv %
\thanks{Jianqing Fan is Frederick L. Moore '18 Professor of Finance, Department of Operations Research and Financial Engineering, Princeton University, Princeton, NJ 08544, USA (e-mail: jqfan@princeton.edu). Jinchi Lv is Assistant Professor of Statistics, Information and Operations Management Department, Marshall School of Business, University of Southern California, Los Angeles, CA 90089, USA (e-mail: jinchilv@marshall.usc.edu).
Fan's research was partially supported by NSF Grants DMS-0704337 and DMS-0714554 and NIH Grant
R01-GM072611. Lv's research was partially supported by NSF Grant DMS-0806030 and 2008 Zumberge Individual Award from USC's James H. Zumberge Faculty Research and Innovation Fund.
We sincerely thank the Co-Editor, Professor Peter Hall, for his kind invitation to write this article. We are also very grateful to the helpful comments of the Co-Editor, Associate Editor and referees that substantially improved the presentation of the paper.}
\medskip\\
Princeton University and University of Southern California\\
} %
}

\maketitle

\begin{abstract}
High dimensional statistical problems arise from diverse fields of scientific research and technological development.   Variable selection plays a pivotal role in contemporary statistical learning and scientific discoveries.  The traditional idea of best subset selection methods, which can be regarded as a specific form of penalized likelihood, is computationally too expensive for many modern statistical applications.  Other forms of penalized likelihood methods have been successfully developed over the last decade to cope with high dimensionality.  They have been widely applied for simultaneously selecting important variables and estimating their effects in high dimensional statistical inference.  In this article, we present a brief account of the recent developments of theory, methods, and implementations for high dimensional variable selection.  What limits of the dimensionality such methods can handle, what the role of penalty functions is, and what the statistical properties are rapidly drive the advances of the field.  The properties of non-concave penalized likelihood and its roles in high dimensional statistical modeling are emphasized.  We also review some recent advances in ultra-high dimensional variable selection, with emphasis on independence screening and two-scale methods.
\end{abstract}

\textit{Short title}: Variable Selection in High Dimensional Feature Space



\newpage

\section{Introduction} \label{Sec1}

High dimensional data analysis has become increasingly frequent and
important in diverse fields of sciences, engineering, and humanities,
ranging from genomics and health sciences to economics, finance and
machine learning.  It characterizes many contemporary problems in
statistics (\cite{HTF09}).   For example, in disease classification using
microarray or proteomics data, tens of thousands of expressions of
molecules or ions are potential predictors; in genowide association studies between
genotypes and phenotypes, hundreds of thousands of SNPs are
potential covariates for phenotypes such as cholesterol levels or
heights. When interactions are considered, the dimensionality grows
quickly. For example, in portfolio allocation among two thousand stocks,
it involves already over two million parameters in the covariance matrix; interactions of molecules in the above examples result in ultra-high dimensionality. To be more precise, throughout the paper ultra-high dimensionality refers to the case where the dimensionality grows at a non-polynomial rate as the sample size increases, and high dimensionality refers to the general case of growing dimensionality. Other examples of high dimensional data include high-resolution images, high-frequency financial data, e-commerce
data, warehouse data, functional, and longitudinal data, among
others. \cite{Donoho00} convincingly demonstrates the
need for developments in high dimensional data analysis, and presents
the curses and blessings of dimensionality. \cite{FL06} give a
comprehensive overview of statistical challenges with high
dimensionality in a broad range of topics, and in particular,
demonstrate that for a host of statistical problems, the model
parameters can be estimated as well as if the best model is known in
advance, as long as the dimensionality is not excessively high. The
challenges that are not present in smaller scale studies have been
reshaping statistical thinking, methodological development, and
theoretical studies.

Statistical accuracy, model interpretability, and computational
complexity are three important pillars of any statistical procedures.
In conventional studies, the number of observations $n$ is much
larger than the number of variables or parameters $p$. In such
cases, none of the three aspects needs to be sacrificed for the
efficiency of others. The traditional methods, however, face
significant challenges when the dimensionality $p$ is comparable to
or larger than the sample size $n$. These challenges include
how to design statistical procedures that are more efficient in
inference; how to derive the asymptotic or nonasymptotic theory; how
to make the estimated models interpretable; and how to make the
statistical procedures computationally efficient and robust.

\begin{figure}
\begin{center}
\includegraphics[height=2 in]
{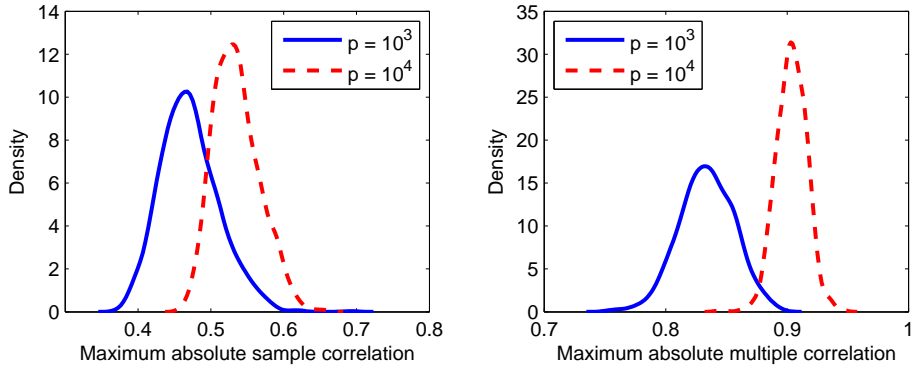}
\begin{singlespace}
\caption{\small Distributions (left panel) of the maximum absolute sample
correlation coefficient $\max_{2 \leq j \leq p} |\mbox{corr}(Z_1,
Z_j)|$, and distributions (right panel) of the maximum absolute multiple correlation coefficient of $Z_1$ with 5 other variables
($\max_{|S| = 5} |\mbox{corr}(Z_1, {\bZ}_S^T \hat{\bbeta}_S)|$, where $\hat{\bbeta}_S$ is the regression coefficient of $Z_1$ regressed on $\bZ_S$, a subset of variables indexed by $S$ and excluding $Z_1$), computed by the stepwise addition algorithm (the actual values are larger than what are presented here), when $n=50$, $p=1000$ (solid curve) and  $p=10000$ (dashed), based on 1000 simulations.
} \label{fig1}
\end{singlespace}
\end{center}%
\end{figure}%

A notorious difficulty of high dimensional model selection comes
from the collinearity among the predictors. The collinearity can
easily be spurious in high dimensional geometry (\cite{FL08}), which
can make us select a wrong model.  Figure~\ref{fig1} shows the
maximum sample correlation and multiple correlation with a given predictor despite that
predictors are generated from independent Gaussian random
variables. As a result, any variable can be well-approximated even by a couple of spurious variables, and can even be replaced by them when the dimensionality is much higher than the sample size. If that variable is a signature predictor and is replaced by spurious variables, we choose wrong variables to associate the covariates with the response and, even worse, the spurious variables can be independent of the response at population level, leading to completely wrong scientific conclusions. Indeed, when the dimensionality $p$ is large, intuition might not be accurate. This is also exemplified by the
data piling problems in high dimensional space observed in
\cite{HMN05}.  Collinearity also gives rise to issues of
over-fitting and model mis-identification.

Noise accumulation in high dimensional prediction has long been recognized in statistics and computer sciences.  Explicit characterization of this is well-known for high dimensional regression problems.  The quantification of the impact of dimensionality on classification was not well understood until \cite{FF08}, who give a simple expression on how dimensionality impacts misclassification rates. \cite{HPG08b} study a similar problem for distanced based-classifiers and observe implicitly the adverse impact of dimensionality. As shown in \cite{FF08}, even for the independence
classification rule described in Section \ref{Sec4.2}, classification using all features can be as bad as a random guess due to noise accumulation in estimating
the population centroids in high dimensional feature space.
Therefore, variable selection is fundamentally important to high dimensional statistical modeling, including regression and classification.

What makes high dimensional statistical inference possible is
the assumption that the regression function lies in a low
dimensional manifold.  In such cases, the $p$-dimensional regression
parameters are assumed to be sparse with many components being zero, where nonzero components indicate the important
variables. With sparsity, variable selection can improve the
estimation accuracy by effectively identifying the subset of
important predictors and can enhance the model interpretability with
parsimonious representation. It can also help reduce the
computational cost when sparsity is very high.

This notion of sparsity is in a narrow sense. It should be
understood more widely in transformed or enlarged feature
spaces. For instance, some prior knowledge may lead us to apply some
grouping or transformation of the input variables (see, e.g.,
\cite{FL08}). Some transformation of the variables may be
appropriate if a significant portion of the pairwise correlations
are high. In some cases, we may want to enlarge the feature space by
adding interactions and higher order terms to reduce the bias of the
model. Sparsity can also be viewed in the context of dimensionality
reduction by introducing a sparse representation, i.e., by reducing the
number of effective parameters in estimation. Examples include the
use of a factor model for high dimensional covariance matrix
estimation in \cite{FFL08}.

Sparsity arises in many scientific endeavors.  In genomic studies,
it is generally believed that only a fraction of molecules are
related to biological outcomes.  For example, in disease
classification, it is commonly believed that only tens of genes
are responsible for a disease.  Selecting tens of genes helps not
only statisticians in constructing a more reliable classification
rule, but also biologists to understand molecular mechanisms.  In
contrast,  popular but naive methods used in microarray data
analysis (\cite{DSB03, ST03, FR06, Efron07}) rely on two-sample
tests to pick important genes, which is truly a marginal correlation
ranking (\cite{FL08}) and can miss important signature genes
(\cite{FSW09b}). The main goals of high dimensional regression and
classification, according to \cite{Bickel08}, are
\begin{itemize}
\item to construct as effective a method as possible to predict future
observations;
\item to gain insight into the relationship between features and
response for scientific purposes, as well as, hopefully, to
construct an improved prediction method.
\end{itemize}
The former appears in problems such as text and document
classification or portfolio optimization, whereas the latter appears
naturally in many genomic studies and other scientific endeavors.

As pointed out in \cite{FL06}, it is helpful to differentiate two
types of statistical endeavors in high dimensional statistical
learning: accuracy of estimated model parameters and accuracy of the
expected loss of the estimated model. The latter property is called
persistence in \cite{GR04} and \cite{Greenshtein06}, and arises
frequently in machine learning problems such as document
classification and computer vision. The former appears in many
other contexts where we want to identify the significant
predictors and characterize the precise contribution of each to the
response variable. Examples include health studies, where the
relative importance of identified risk factors needs to be assessed
for prognosis. Many of the existing results in the literature have
been concerned with the study of consistency of high dimensional variable
selection methods, rather than characterizing the asymptotic
distributions of the estimated model parameters. However, consistency and
persistence results are inadequate for understanding uncertainty
in parameter estimation.

High dimensional variable selection encompasses a majority of
frontiers where statistics advances rapidly today.  There has been
an evolving literature in the last decade devoted to understanding the performance of
various variable selection techniques. The main theoretical questions include determining the limits of the dimensionality that such methods can handle and how to characterize the optimality of variable selection procedures. The answers to the first question for many
existing methods were largely unknown until recently.  To a large
extent, the second question still remains open for many procedures.
In the Gaussian linear regression model, the case of orthonormal
design reduces to the problem of Gaussian mean estimation, as do the
wavelet settings where the design matrices are orthogonal. In such
cases, the risks of various shrinkage estimators and their optimality
have been extensively studied. See, e.g., \cite{DJ94} and
\cite{AF01}.

In this article we address the issues of variable selection for
high dimensional statistical modeling in the unified framework of
penalized likelihood estimation. It has been widely used in
statistical inferences and machine learning, and is basically a
moderate scale learning technique.  We also give an overview on the techniques for ultrahigh dimensional screening.  Combined iteratively with large scale screening, it can handle problems of ultra-high dimensionality
(\cite{FSW09b}).  This will be reviewed as well.

The rest of the article is organized as follows. In Section \ref{Sec2}, we discuss the connections of penalized likelihood to classical model selection methods. Section \ref{Sec3} details the methods and implementation of penalized likelihood estimation. We review some recent advances in ultra-high dimensional variable selection in Section \ref{Sec4}. In Section \ref{Sec5}, we survey the sampling properties of penalized least squares. Section \ref{Sec6} presents the classical oracle property of penalized least squares and penalized likelihood methods in ultra-high dimensional space. We conclude the article with some additional remarks in Section \ref{Sec7}.

\section{Classical model selection} \label{Sec2}

Suppose that the available data are $(\bx_i\t, y_i)_{i = 1}^n$, where
$y_i$ is the $i$-th observation of the response variable and $\bx_i$ is its associated
$p$-dimensional covariates vector.  They are usually assumed to be a random
sample from the population $(\bX\t, Y)$, where the conditional mean of $Y$ given $\bX$
depends on the linear predictor $\bbeta^T\bX$ with $\bbeta = (\beta_1, \cdots, \beta_p)\t$.  In sparse modeling, it is frequently assumed that most regression coefficients $\beta_j$
are zero. Variable selection aims to identify all important
variables whose regression coefficients do not vanish and to provide
effective estimates of those coefficients.

More generally, assume that the data are generated from the true density function $f_{\btheta_0}$ with parameter
vector $\btheta_0 = (\theta_1, \cdots, \theta_d)\t$. Often, we are
uncertain about the true density, but more certain about a larger family of models $f_{\btheta_1}$ in which $\btheta_0$ is a (nonvanishing) subvector of the $p$-dimensional parameter vector $\btheta_1$.
The problems of how to estimate the dimension of the model and
compare models of different dimensions naturally arise in many
statistical applications, including time series modeling. These are referred to as model selection in the literature.

\cite{Akaike73, Akaike74} proposes to choose a model that minimizes
the Kullback-Leibler (KL) divergence of the fitted model from the true model.
\cite{Akaike73} considers the maximum likelihood estimator (MLE)
$\hbtheta = (\hat \theta_1, \cdots, \hat \theta_p)\t$ of the parameter vector $\btheta$ and shows that, up to
an additive constant, the estimated KL divergence can be
asymptotically expanded as
$$
    -\ell_n(\hbtheta) + \lambda \dim(\hbtheta) =   -\ell_n(\hbtheta) + \lambda \sum_{j=1}^p
          I(\hat \theta_j \not = 0),
$$
where $\ell_n(\btheta)$ is the log-likelihood function,
$\dim(\btheta)$ denotes the dimension of the  model, and $\lambda =
1$. This leads to the AIC. \cite{Schwartz78} takes a Bayesian approach with prior
distributions that have nonzero prior probabilities on some lower
dimensional subspaces and proposes the BIC with $\lambda = (\log
n)/2$ for model selection. Recently, \cite{LL08}
gave a KL divergence interpretation of Bayesian model selection and
derive generalizations of AIC and BIC when the model may be misspecified.

The work of AIC and BIC suggests a unified approach to model
selection: choose a parameter vector $\btheta$ that maximizes the
penalized likelihood
\begin{equation} \label{e01}
\ell_n(\btheta) - \lambda \|\btheta\|_0,
\end{equation}
where the $L_0$-norm of $\btheta$ counts the number of non-vanishing
components in $\btheta$ and $\lambda \geq 0$ is a regularization parameter.  Given $\|\btheta\|_0 = m$, the solution to
(\ref{e01}) is the subset with the largest maximum
likelihood among all subsets of size $m$.  The model size $m$ is then
chosen to maximize (\ref{e01}) among $p$ best subsets of sizes $m$ ($1
\leq m \leq p$).  Clearly, the computation of the penalized $L_0$
problem is a combinational problem with NP-complexity.

When the normal likelihood is used, (\ref{e01}) becomes penalized
least squares. Many traditional methods can be regarded as
penalized likelihood methods with different choices of $\lambda$. Let
$\text{RSS}_d$ be the residual sum of squares of the best subset
with $d$ variables.   Then $C_p = \text{RSS}_d/s^2 + 2 d - n$ in
\cite{Mallows73} corresponds to $\lambda = 1$, where $s^2$ is the mean squared error of the full model.  The adjusted $R^2$
given by
\[
   R^2_{\text{adj}} = 1 - \frac{n - 1}{n - d} \frac{\text{RSS}_d}{\text{SST}}
\]
also amounts to a penalized-$L_0$ problem, where SST is the total sum of squares. Clearly maximizing
$R^2_{\text{adj}}$ is equivalent to minimizing $\log(\text{RSS}_d/(n -
d))$. By $\text{RSS}_d/n \approx \sigma^2$ (the error variance), we have
\[
   n \log \frac{\text{RSS}_d}{n - d}
   \approx \text{RSS}_d /\sigma^2 + d + n (\log \sigma^2 -1).
\]
This shows that the adjusted $R^2$ method is approximately equivalent to PMLE with $\lambda = 1/2$.  Other examples include the
generalized cross-validation (GCV) given by $\text{RSS}_d/(1 -
d/n)^2$, cross-validation (CV), and RIC in \cite{FG94}.
See \cite{BL06} for more discussions of regularization in
statistics.

\section{Penalized likelihood} \label{Sec3}

As demonstrated above, $L_0$ regularization arises naturally in
many classical model selection methods. It gives a nice interpretation
of best subset selection and admits nice sampling properties
(\cite{BBM99}).  However, the computation is infeasible in
high dimensional statistical endeavors.  Other penalty functions
should be used.  This results in a generalized form
\begin{equation} \label{e02}
n^{-1} \ell_n(\bbeta) - \sum_{j=1}^p p_\lambda (|\beta_j|),
\end{equation}
where $\ell_n(\bbeta)$ is the log-likelihood function and $p_\lambda(\cdot)$ is a penalty function indexed by the regularization parameter $\lambda \geq 0$. By maximizing the
penalized likelihood (\ref{e02}), we hope to simultaneously select
variables and estimate their associated regression coefficients.  In
other words, those variables whose regression coefficients are
estimated as zero are automatically deleted.

A natural generalization of penalized $L_0$-regression is penalized
$L_q$-regression, called bridge regression in \cite{FF93}, in
which $p_\lambda(|\theta|) = \lambda |\theta|^q$ for $0 < q \leq
2$.   This bridges the best subset section (penalized $L_0$) and
ridge regression (penalized $L_2$), including the $L_1$-penalty as a
specific case. The non-negative garrote is introduced in \cite{Breiman95} for shrinkage estimation and variable selection. Penalized $L_1$-regression is called the LASSO by
\cite{Tibshirani96} in the ordinary regression setting, and is now
collectively referred to as penalized $L_1$-likelihood.  Clearly,
penalized $L_0$-regression possesses the variable selection feature,
whereas penalized $L_2$-regression does not.  What kind of
penalty functions are good for model selection?

\cite{FL01} advocate penalty functions that
give estimators with three properties:
\begin{itemize}
\item[1)] \textit{Sparsity}: The resulting estimator automatically
sets small estimated coefficients to zero to accomplish variable
selection and reduce model complexity.

\item[2)] \textit{Unbiasedness}: The resulting estimator
is nearly unbiased, especially when the true coefficient
$\beta_j$ is large, to reduce model bias.

\item[3)] \textit{Continuity}: The resulting estimator is continuous
in the data to reduce instability in model
prediction (\cite{Breiman96}).
\end{itemize}
They require the penalty function $p_\lambda(|\theta|)$ to be nondecreasing in $|\theta|$, and
provide insights into these properties.  We first consider the penalized least squares in a canonical form.

\subsection{Canonical regression model} \label{Sec3.1}

Consider the linear regression model
\begin{equation} \label{e03}
\by = \bX \bbeta + \bveps,
\end{equation}
where $\bX = (\bx_1, \cdots, \bx_n)\t$, $\by = (y_1, \cdots,
y_n)\t$, and $\bveps$ is an $n$-dimensional noise vector. If
$\bveps \sim N(\bzero, \sigma^2 I_n)$, then the penalized likelihood (\ref{e02}) is
equivalent, up to an affine transformation of the log-likelihood, to the penalized least squares (PLS) problem
\begin{equation} \label{e04}
\min_{\bbeta \in \mathbf{R}^p} \left\{ \frac{1}{2 n}  \|\by - \bX \bbeta\|^2 + \sum_{j=1}^p p_\lambda (|\beta_j|)\right\},
\end{equation}
where $\|\cdot\|$ denotes the $L_2$-norm. Of course, the penalized least squares continues to be applicable
even when the noise does not follow a normal distribution.

For the canonical linear model in which the design matrix multiplied by $n^{-1/2}$ is
orthonormal (i.e., $\bX\t \bX = n I_p$), (\ref{e04}) reduces to the minimization of
\begin{equation} \label{e05}
 \frac{1}{2 n} \|\by - \bX \hbbeta\|^2 + \|\hbbeta - \bbeta\|^2 + \sum_{j=1}^p p_\lambda (|\beta_j|),
\end{equation}
where $\hbbeta = n^{-1} \bX\t \by$ is the ordinary least squares estimate.
Minimizing (\ref{e05}) becomes a componentwise regression
problem.  This leads to considering the univariate PLS problem
\begin{equation} \label{e06}
\hat{\theta}(z) = \arg\min_{\theta \in \mathbf{R}}
\left\{\frac{1}{2} (z - \theta)^2 + p_\lambda(|\theta|)\right\}.
\end{equation}
\cite{AF01} show that the PLS estimator $\hat{\theta}(z)$
possesses the properties:
\begin{itemize}
\item[1)] \textit{sparsity} if $\min_{t \geq 0} \{t + p_\lambda'(t)\} > 0$;

\item[2)] \textit{approximate unbiasedness} if $p_\lambda'(t) = 0$ for large $t$;

\item[3)] \textit{continuity} if and only if $\arg\min_{t \geq 0} \{t + p_\lambda'(t)\} = 0$,
\end{itemize}
where $p_\lambda(t)$ is nondecreasing and continuously differentiable on $[0, \infty)$, the function $-t - p_\lambda'(t)$ is strictly unimodal on $(0, \infty)$, and $p_\lambda'(t)$ means $p_\lambda'(0+)$ when $t = 0$ for notational simplicity. In general for the penalty function, the singularity at the origin (i.e., $p_\lambda'(0+) > 0$) is needed for generating sparsity in variable selection and the concavity is needed to reduce the estimation bias.

\subsection{Penalty function} \label{Sec3.2}

\begin{figure}
\begin{center}
\includegraphics[height=2 in]
{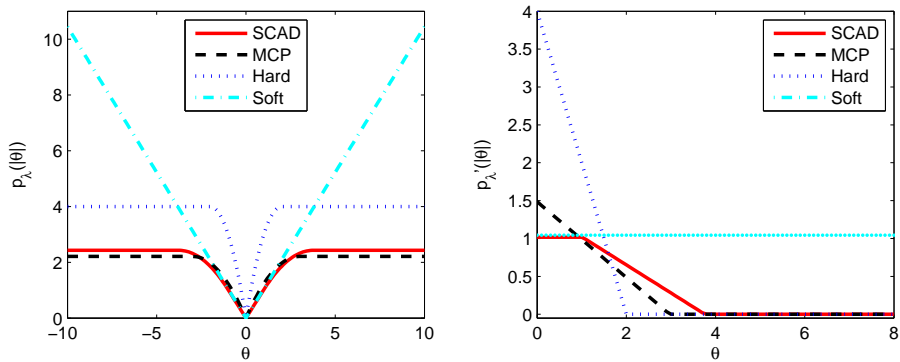}
\begin{singlespace}
\caption{\small Some commonly used penalty functions (left panel)
and their derivatives (right panel).  They correspond to the risk
functions shown in the right panel of Figure~\ref{fig3}.  More
precisely, $\lambda = 2$ for hard thresholding penalty, $\lambda =
1.04$ for $L_1$-penalty, $\lambda = 1.02$ for SCAD with $a=3.7$, and
$\lambda=1.49$ for MCP with $a=2$.} \label{fig2}
\end{singlespace}
\end{center}%
\end{figure}%

It is known that the convex $L_q$ penalty with $q > 1$ does not
satisfy the sparsity condition, whereas the convex $L_1$ penalty
does not satisfy the unbiasedness condition, and the concave $L_q$
penalty with $0 \leq q < 1$ does not satisfy the continuity
condition. In other words, none of the $L_q$ penalties satisfies all three conditions simultaneously.  For this reason,
\cite{Fan97} and \cite{FL01} introduce the smoothly clipped absolute deviation
(SCAD), whose derivative is given by
\begin{equation} \label{e07}
p_\lambda'(t) = \lambda \left\{I\left(t \leq \lambda\right)
+ \frac{\left(a \lambda - t\right)_+}{\left(a - 1\right) \lambda}
I\left(t > \lambda\right)\right\} \quad \text{for some } a > 2,
\end{equation}
where $p_\lambda(0) = 0$ and, often, $a = 3.7$ is used (suggested by a Bayesian
argument).  It satisfies the aforementioned three properties.  A
penalty of similar spirit is the minimax concave penalty (MCP) in
\cite{Zhang09}, whose derivative is given by
\begin{equation}  \label{e08}
p_\lambda'(t) =\left(a \lambda - t\right)_+/a.
\end{equation}
Clearly SCAD takes off at the origin as the $L_1$ penalty and then
gradually levels off, and MCP translates the flat part of the derivative of
SCAD to the origin.  When
\begin{equation}  \label{e09}
   p_\lambda (t) = \lambda^2 - (\lambda - t)_+^2,
\end{equation}
\cite{Antoniadis96} shows that the solution is the hard-thresholding estimator $\hat{\theta}_H(z) = z I(|z| > \lambda)$.  A family of concave penalties that bridge the
$L_0$ and $L_1$ penalties is studied by \cite{LF09} for model
selection and sparse recovery.  A linear combination of $L_1$ and $L_2$ penalties is called an elastic net by \cite{ZH05}, which encourages some grouping effects. Figure~\ref{fig2} depicts some of those commonly used penalty functions.

We now look at the PLS estimator $\hat \theta(z)$ in (\ref{e06}) for a
few penalties. Each increasing penalty function gives a shrinkage
rule: $|\hat \theta(z)| \leq |z|$ and $\hat \theta(z) = \sgn(z) |\hat \theta(z)|$ (\cite{AF01}). The entropy penalty ($L_0$ penalty)
and the hard thresholding penalty yield the hard thresholding rule
(\cite{DJ94}), while the $L_1$ penalty gives the soft thresholding
rule (\cite{Bickel83, DJ94}). The SCAD and MCP give rise to
analytical solutions to (\ref{e06}), each of which is a linear spline in $z$ (\cite{Fan97}).

\begin{figure}
\begin{center}
\includegraphics[height=2in]
{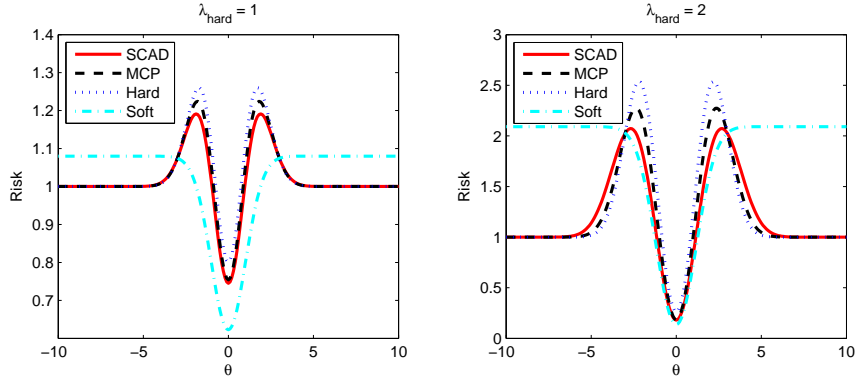}
\begin{singlespace}
\caption{\small The risk functions for penalized least squares under
the Gaussian model for the hard-thresholding penalty, $L_1$-penalty,
SCAD ($a = 3.7$), and MCP ($a = 2$). The left panel corresponds to $\lambda =
1$ and the right panel corresponds to $\lambda =2$ for the
hard-thresholding estimator, and the rest of parameters are chosen so
that their risks are the same at the point $\theta = 3$.}
\label{fig3}
\end{singlespace}
\end{center}%
\end{figure}%

How do those thresholded-shrinkage estimators perform?  To compare
them, we compute their risks in the fundamental model in which $Z
\sim N(\theta, 1)$.  Let $R(\theta) = E (\hat{\theta}(Z) -
\theta)^2$. Figure~\ref{fig3} shows the risk functions $R(\theta)$
for some commonly used penalty functions.  To make them comparable, we chose
$\lambda = 1$ and $2$ for the hard thresholding penalty, and for other penalty functions the values of $\lambda$ were chosen to make their risks at $\theta = 3$ the same.
Clearly the penalized likelihood estimators improve the ordinary
least squares estimator $Z$ in the region where $\theta$ is near
zero, and have the same risk as the ordinary least squares estimator
when $\theta$ is far away from zero (e.g., 4 standard
deviations away), except the LASSO estimator. When $\theta$ is large, the LASSO estimator has a bias approximately of size $\lambda$, and this causes higher
risk as shown in Figure~\ref{fig3}. When
$\lambda_\text{hard} = 2$, the LASSO estimator has higher risk than the SCAD
estimator, except in a small region.  The bias of the LASSO estimator
makes LASSO prefer a smaller $\lambda$. For $\lambda_\text{hard} = 1$, the
advantage of the LASSO estimator around zero is more pronounced.  As a
result in model selection, when $\lambda$ is automatically
selected by a data-driven rule to compensate the bias problem, the
LASSO estimator has to choose a smaller $\lambda$ in order to have a
desired mean squared error. Yet, a smaller value of $\lambda$ results
in a complex model.  This explains why the LASSO estimator tends
to have many false positive variables in the selected model.

\subsection{Computation and implementation} \label{Sec3.3}

It is challenging to solve the penalized likelihood problem (\ref{e02})
when the penalty function $p_\lambda$ is nonconvex.   Nevertheless, \cite{FL09} are able to give the conditions under which the penalized likelihood estimator exists and is unique; see also \cite{KK09} for the results of penalized least squares with SCAD penalty. When the
$L_1$-penalty is used, the objective function (\ref{e02}) is concave and hence convex optimization algorithms can
be applied.  We show in this section that the penalized likelihood
(\ref{e02}) can be solved by a sequence of reweighted penalized
$L_1$-regression problems via local linear approximation
(\cite{ZL08}).

In the absence of other available algorithms at that time,  \cite{FL01}
propose a unified and effective local quadratic approximation (LQA)
algorithm for optimizing nonconcave penalized likelihood. Their
idea is to locally approximate the objective function by a quadratic
function. Specifically, for a given initial
value $\bbeta^* = (\beta^*_1, \cdots, \beta^*_p)\t$, the penalty function $p_\lambda$ can be
locally approximated by a quadratic function as
\begin{equation} \label{e10}
  p_\lambda(|\beta_j|) \approx p_\lambda(|\beta^*_j|)
  + \frac{1}{2} \frac{p'_\lambda(|\beta^*_j|)}{|\beta^*_j|} [\beta_j^2 - (\beta^*_j)^2] \quad \text{ for } \beta_j \approx \beta^*_j.
\end{equation}
With this and a LQA to the log-likelihood, the penalized likelihood (\ref{e02}) becomes a least squares problem that admits a closed-form solution. To avoid numerical instability, it
sets the estimated coefficient $\hbeta_j = 0$ if $\beta^*_j$ is very
close to 0, which amounts to deleting the $j$-th covariate from the
final model.
Clearly the value $0$ is an absorbing state of LQA in the sense that once a coefficient is set to zero, it remains zero in subsequent iterations.

The convergence property of the LQA was studied in \cite{HL05}, who
show that LQA plays the same role as the E-step in the EM
algorithm in \cite{Dempster77}. Therefore LQA has similar behavior
to EM. Although the EM requires a full iteration for maximization
after each E-step, the LQA updates the quadratic approximation at
each step during the course of iteration, which speeds up the
convergence of the algorithm. The convergence rate of LQA is
quadratic, which is the same as that of the modified EM algorithm in
\cite{Lange95}.

\begin{figure}
\begin{center}
\includegraphics[height=2in]
{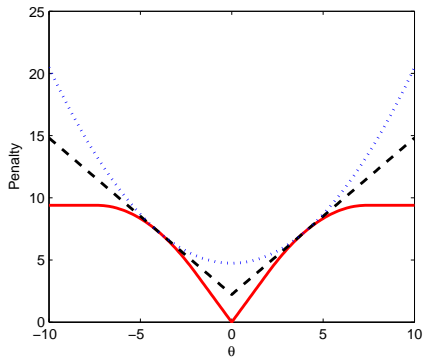}
\begin{singlespace}
\caption{\small The local linear (dashed) and local quadratic (dotted) approximations
to the SCAD function (solid) with $\lambda = 2$ and $a = 3.7$ at a given point $|\theta| = 4$.} \label{fig4}
\end{singlespace}
\end{center}%
\end{figure}%

A better approximation can be achieved by using the local linear
approximation (LLA):
\begin{equation} \label{e11}
p_\lambda(|\beta_j|) \approx p_\lambda(|\beta^*_j|) +
p'_\lambda(|\beta^*_j|)(|\beta_j| - |\beta^*_j|) \quad \text{ for }
\beta_j \approx \beta^*_j,
\end{equation}
as in \cite{ZL08}. See Figure~\ref{fig4} for an illustration of the
local linear and local quadratic approximations to the SCAD function.
Clearly, both LLA and LQA are convex majorants of concave penalty function $p_\lambda(\cdot)$ on
$[0, \infty)$, but LLA is a better approximation since it is the
minimum (tightest) convex majorant of the concave function on $[0, \infty)$.
With LLA, the penalized likelihood (\ref{e02}) becomes
\begin{equation}   \label{e12}
n^{-1} \ell_n(\bbeta) - \sum_{j = 1}^p w_j |\beta_j|,
\end{equation}
where the weights are $w_j = p'_\lambda(|\beta^*_j|)$.  Problem
(\ref{e12}) is a concave optimization problem if the log-likelihood
function is concave.  Different penalty functions give different
weighting schemes, and LASSO gives a constant weighting scheme.  In
this sense, the nonconcave penalized likelihood is an
iteratively reweighted penalized $L_1$ regression.  The weight
function is chosen adaptively to reduce the biases due to
penalization.  For example, for SCAD and MCP, when the estimate of a
particular component is large so that it has high confidence to be
non-vanishing, the component does not receive any penalty in
(\ref{e12}), as desired.

\cite{Zou06} proposes the weighting scheme $w_j =
|\beta^*_j|^{-\gamma}$ for some $\gamma > 0$ and calls the
resulting procedure adaptive LASSO.  This weight reduces the penalty
when the previous estimate is large.  However, the penalty at zero
is infinite.  When the procedure is applied iteratively, zero
becomes an absorbing state. On the
other hand, the penalty functions such as SCAD and MCP do not have
this undesired property. For example, if the initial estimate is
zero, then $w_j = \lambda$ and the resulting estimate is the LASSO estimate.

\cite{FL01}, \cite{Zou06}, and \cite{ZL08} all suggest a
consistent estimate such as the un-penalized MLE.  This implicitly
assumes that $p \ll n$.  For dimensionality $p$ that is larger than
sample size $n$, the above method is not applicable.  \cite{FL08} recommend
using $\beta^*_j = 0$, which is equivalent to using the LASSO estimate as
the initial estimate.  Another possible initial value is to use a stepwise addition fit or componentwise regression. They put forward the recommendation that only
a few iterations are needed, which is in line with \cite{ZL08}.

Before we close this section, we remark that with the LLA and LQA, the resulting sequence of target values is always nondecreasing, which is a specific feature of minorization-maximization (MM) algorithms (\cite{HL00}).  Let $p_\lambda(\bbeta) = \sum_{j=1}^p p_\lambda(|\beta_j|)$. Suppose that at the $k$-th iteration, $p_\lambda(\bbeta)$ is approximated by $q_{\lambda}(\bbeta)$ such that
\begin{equation} \label{e13}
    p_\lambda (\bbeta) \leq q_{\lambda} (\bbeta) \quad \text{ and } \quad
    p_\lambda (\bbeta^{(k)}) = q_{\lambda} (\bbeta^{(k)}),
\end{equation}
where $\bbeta^{(k)}$ is the estimate at the $k$-th iteration.  Let $\bbeta^{(k+1)}$ maximize the approximated penalized likelihood $n^{-1} \ell_n(\bbeta) - q_{\lambda} (\bbeta)$.  Then we have
\begin{eqnarray*}
n^{-1} \ell_n(\bbeta^{(k+1)}) - p_{\lambda} (\bbeta^{(k+1)}) & \geq  &
n^{-1} \ell_n(\bbeta^{(k+1)}) - q_{\lambda} (\bbeta^{(k+1)})\\
& \geq &
n^{-1} \ell_n(\bbeta^{(k)}) - q_{\lambda} (\bbeta^{(k)})\\
& =  & n^{-1} \ell_n(\bbeta^{(k)}) - p_{\lambda} (\bbeta^{(k)}).
\end{eqnarray*}
Thus, the target values are non-decreasing.  Clearly, the LLA and LQA are two specific cases of the MM algorithms, satisfying condition (\ref{e13}); see Figure~\ref{fig4}.  Therefore, the sequence of target function values is non-decreasing and thus converges provided it is bounded.  The critical point is the global maximizer under the conditions in \cite{FL09}.

\subsection{LARS and other algorithms} \label{Sec3.4}

As demonstrated in the previous section, the penalized least squares problem (\ref{e04}) with an $L_1$ penalty is fundamental to the computation of penalized likelihood estimation.  There are several additional powerful algorithms for such an endeavor.  \cite{OPT00} cast such a problem as a quadratic programming problem.  \cite{EHJT04} propose a fast and efficient
least angle regression (LARS) algorithm for variable selection, a
simple modification of which produces the entire LASSO solution path
$\{\widehat{\bbeta}(\lambda): \lambda>0\}$ that optimizes (\ref{e04}).  The computation is based on the fact that the LASSO solution path is piecewise linear in $\lambda$.  See \cite{RZ07} for a more general account of the conditions under which the solution to the penalized likelihood  (\ref{e02}) is piecewise linear.  The LARS algorithm starts from a large value of $\lambda$ which selects only one covariate that has the greatest correlation with the response variable and decreases the $\lambda$ value until the second variable is selected, at which the selected variables have the same correlation (in magnitude) with the current working residual as the first one, and so on. See  \cite{EHJT04} for details.

The idea of the LARS algorithm can be expanded to compute the solution paths of penalized least squares (\ref{e04}).  \cite{Zhang09} introduces the PLUS algorithm for efficiently computing a solution path of (\ref{e04}) when the penalty function $p_\lambda(\cdot)$ is a quadratic spline such as the SCAD and MCP.  In addition, \cite{Zhang09} also shows that the solution path $\widehat{\bbeta}(\lambda)$ is piecewise linear in $\lambda$, and the proposed solution path has desired statistical properties.

For the penalized least squares problem (\ref{e04}), \cite{Fu98},
\cite{DDD04}, and \cite{WL08} propose a coordinate descent algorithm,
which iteratively optimizes (\ref{e04}) one component at a time.
This algorithm can also be applied to optimize the group LASSO
(\cite{AF01, YL06}) as shown in \cite{MvB08}, penalized precision matrix estimation
(\cite{FHT07}), and penalized likelihood (\ref{e02}) (\cite{FL09, ZL09}).

More specifically, \cite{FL09} employ a path-following coordinate optimization algorithm, called the iterative coordinate ascent (ICA) algorithm, for maximizing the nonconcave penalized likelihood. It successively maximizes the penalized likelihood (\ref{e02}) for regularization parameters $\lambda$ in decreasing order. A similar idea is also studied in \cite{ZL09}, who introduce the ICM algorithm. The coordinate optimization algorithm uses the Gauss-Seidel method, i.e., maximizing one coordinate at a time with successive displacements. Specifically, for each coordinate within each iteration, it uses the second order approximation of $\ell_n(\bbeta)$ at the $p$-vector from the previous step along that coordinate and maximizes the univariate penalized quadratic approximation
\begin{equation} \label{e50}
\max_{\theta \in \mathbf{R}} \left\{-\frac{1}{2} (z - \theta)^2 - \Lambda p_\lambda(|\theta|)\right\},
\end{equation}
where $\Lambda > 0$. It updates each coordinate if the maximizer of the corresponding univariate penalized quadratic approximation makes the penalized likelihood (\ref{e02}) strictly increase. Therefore, the ICA algorithm enjoys the ascent property that the resulting sequence of values of the penalized likelihood is increasing for a fixed $\lambda$. Compared to other algorithms, the coordinate optimization algorithm is especially appealing for large scale problems with both $n$ and $p$ large, thanks to its low computational complexity. It is fast to implement when the univariate problem (\ref{e50}) admits a closed-form solution. This is the case for many commonly used penalty functions such as SCAD and MCP. In practical implementation, we pick a sufficiently large $\lambda_{\max}$ such that the maximizer of the penalized likelihood (\ref{e02}) with $\lambda = \lambda_{\max}$ is $\bzero$, and a decreasing sequence of regularization parameters. The studies in \cite{FL09} show that the coordinate optimization works equally well and efficiently for producing the entire solution paths for concave penalties.

The LLA algorithm for computing penalized likelihood is now available in R at
\begin{center}
http://cran.r-project.org/web/packages/SIS/index.html
\end{center}
as a function in the SIS package.  So is the PLUS algorithm for computing the penalized least squares estimator with SCAD and MC+ penalties.  The Matlab codes are also available for the ICA algorithm for computing the solution path of the penalized likelihood estimator and for computing SIS upon request.

\subsection{Composite quasi-likelihood} \label{Sec3.5}

The function $\ell_n(\bbeta)$ in (\ref{e02}) does not have to be the true likelihood.  It can be a quasi-likelihood or a loss function (\cite{FSW09b}).  In most statistical applications, it is of the form
\begin{equation}   \label{e14}
 n^{-1} \sum_{i=1}^n Q(\bx_i\t\bbeta, y_i) - \sum_{j=1}^p p_\lambda(|\beta_j|),
\end{equation}
where $Q(\bx_i\t\bbeta, y_i)$ is the conditional quasi-likelihood of $Y_i$ given $\bX_i$.  It can also be the loss function of using $\bx_i\t\bbeta$ to predict $y_i$.  In this case, the penalized quasi-likelihood (\ref{e14}) is written as the minimization of
\begin{equation}   \label{e15}
 n^{-1} \sum_{i=1}^n L(\bx_i\t\bbeta, y_i) + \sum_{j=1}^p p_\lambda(|\beta_j|),
\end{equation}
where $L$ is a loss function.  For example, the loss function can be a robust loss:  $L(x, y) = |y-x|$. How should we choose a quasi-likelihood to enhance the efficiency of procedure when the error distribution possibly deviates from normal?

To illustrate the idea, consider the linear model (\ref{e03}).  As long as the error distribution of $\varepsilon$ is homoscedastic, $\bx_i\t\bbeta$ is, up to an additive constant, the conditional $\tau$ quantile of $y_i$ given $\bx_i$.  Therefore, $\bbeta$ can be estimated by the quantile regression
$$
     \sum_{i=1}^n \rho_\tau (y_i - b_\tau - \bx_i\t \bbeta),
$$
where $\rho_\tau(x) = \tau x_+ + (1-\tau) x_-$ (\cite{KB78}).  \cite{Koenker84} proposes solving the weighted composite quantile regression by using different quantiles to improve the efficiency, namely, minimizing with respect to $b_1, \cdots, b_K$ and $\bbeta$,
\begin{equation} \label{e16}
   \sum_{k=1}^K w_k \sum_{i=1}^n \rho_{\tau_k}(y_i - b_k - \bx_i\t \bbeta),
\end{equation}
where $\{\tau_k\}$ is a given sequence of quantiles and $\{w_k\}$ is a given sequence of weights.  \cite{ZY08} propose the penalized composite quantile with equal weights to improve the efficiency of the penalized least squares.

Recently, \cite{BFW09} proposed the more general composite quasi-likelihood
\begin{equation}   \label{e17}
 \sum_{k=1}^K w_k \sum_{i=1}^n L_k(\bx_i\t\bbeta, y_i) + \sum_{j=1}^p p_\lambda(|\beta_j|).
\end{equation}
They derive the asymptotic normality of the estimator and choose the weight function to optimize the asymptotic variance.  In this view, it always performs better than a single quasi-likelihood function.  In particular, they study in detail the relative efficiency of the composite $L_1$-$L_2$ loss and optimal composite quantile loss with the least squares estimator.

Note that the composite likelihood (\ref{e17}) can be regarded as an approximation to the log-likelihood function via
$$
   \log f(y|\bx) = \log f(y|\bx\t\bbeta) \approx -\sum_{k=1}^K w_k L_k(\bx\t\bbeta, y)
$$
with $\sum_{k=1}^K w_k = 1$.  Hence, $w_k$ can also be chosen to minimize (\ref{e17}) directly.  If the convexity of the composite likelihood is enforced, we need to impose the additional constraint that all weights are non-negative.

\subsection{Choice of penalty parameters} \label{Sec3.6}

The choice of penalty parameters is of paramount importance in penalized likelihood estimation.  When $\lambda = 0$, all variables are selected and the model is even unidentifiable when $p > n$.  When $\lambda = \infty$, if the penalty satisfies $\lim_{\lambda \to \infty} p_\lambda (|\theta|) = \infty$ for $\theta \not = 0$, then none of the variables is selected.  The interesting cases lie between these two extreme choices.

The above discussion clearly indicates that $\lambda$ governs the complexity of the selected model.  A large value of $\lambda$ tends to choose a simple model, whereas a small value of $\lambda$ inclines to a complex model.  The estimation using a larger value of $\lambda$ tends to have smaller variance, whereas the estimation using a smaller value of $\lambda$ inclines to smaller modeling biases.  The trade-off between the biases and variances yields an optimal choice of $\lambda$.  This is frequently done by using a multi-fold cross-validation.

There are relatively few studies on the choice of penalty parameters.  In \cite{WLT07}, it is shown that the model selected by generalized cross-validation using the SCAD penalty contains all important variables, but with nonzero probability includes some unimportant variables, and that the model selected by using BIC achieves the model selection consistency and an oracle property. It is worth to point out that missing some true predictor causes model misspecification, as does misspecifying the family of distributions. A semi-Bayesian information criterion (SIC) is proposed by \cite{LL08} to address this issue for model selection.

\section{Ultra-high dimensional variable selection} \label{Sec4}

Variable selection in ultra-high dimensional feature space has
become increasingly important in statistics, and calls for new or
extended statistical methodologies and theory. For example, in
disease classification using microarray gene expression data, the
number of arrays is usually on the order of tens while the number of
gene expression profiles is on the order of tens of thousands; in
the study of protein-protein interactions, the number of features
can be on the order of millions while the sample size $n$ can be
on the order of thousands (see, e.g., \cite{THNC03} and
\cite{FR06}); the same order of magnitude occurs in genetic association studies between genotypes and phenotypes. In such problems,
it is important to identify significant features (e.g., SNPs) contributing to the response and
reliably predict certain clinical prognosis (e.g., survival time and cholesterol level). As mentioned in the introduction, three important
issues arise in such high dimensional statistical endeavors:
computational cost, statistical accuracy, and model
interpretability. Existing variable selection techniques can become
computationally intensive in ultra-high dimensions.

A natural idea is to reduce the dimensionality $p$ from a large or
huge scale (say, $\log p = O(n^a)$ for some $a > 0$) to a relatively
large scale $d$ (e.g., $O(n^b)$ for some $b > 0$) by a fast, reliable, and efficient method, so that well-developed variable selection techniques can be applied to
the reduced feature space. This provides a powerful tool for
variable selection in ultra-high dimensional feature space. It
addresses the aforementioned three issues when the variable
screening procedures are capable of retaining all the important
variables with asymptotic probability one, the sure screening
property introduced in \cite{FL08}.

\begin{figure}
\begin{center}
\includegraphics[height=1in]
{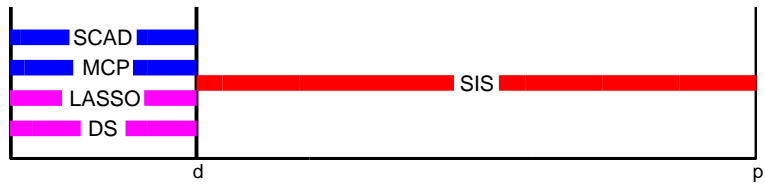}
\begin{singlespace}
\caption{\small Illustration of ultra-high dimensional variable selection scheme.  A large scale screening is first used to screen out unimportant variables and then a moderate-scale searching is applied to further select important variables.  At both steps, one can choose a favorite method. } \label{fig5}
\end{singlespace}
\end{center}%
\end{figure}%

The above discussion suggests already a two-scale method for ultra-high dimensional variable selection problems:  a crude large scale screening followed by a moderate scale selection.   The idea is explicitly suggested by \cite{FL08} and is illustrated by the schematic diagram in Figure~\ref{fig5}.   One can choose one of many popular screening techniques, as long as it possesses the sure screening property.  In the same vein, one can also select a preferred tool for moderate scale selection.  The large-scale screening and moderate-scale selection can be iteratively applied, resulting in iterative sure independence screening (ISIS) (\cite{FL08}).  Its amelioration and extensions are given in \cite{FSW09b}, who also develop R and Matlab codes to facilitate the implementation in generalized linear models (\cite{MN89}).

\subsection{Sure independence screening} \label{Sec4.1}

Independence screening refers to ranking features according to marginal utility, namely, each feature is used independently as a predictor to decide its usefulness for predicting the response.  Sure independence screening (SIS) was introduced by \cite{FL08} to reduce the computation in ultra-high dimensional variable selection:  all important features are in the selected model with probability tending to 1 (\cite{FL08}).  An example of independence learning is the correlation ranking proposed in \cite{FL08} that ranks features according to the magnitude of its sample correlation with the response variable.  More precisely, let $\bomega = (\omega_1, \cdots, \omega_p)\t = \bX\t \by$ be a $p$-vector obtained by componentwise regression, where we assume that each column of the $n \times p$ design matrix $\bX$ has been standardized with mean zero and variance one.  For any given $d_n$, take the selected submodel to be
\begin{equation} \label{e18}
\widehat{\mathcal{M}}_d = \{1 \leq j \leq p: |\omega_j| \text{ is among
the first $d_n$ largest of all}\}.
\end{equation}
This reduces the full model of size $p \gg n$ to a submodel with size $d_n$, which can be less than $n$.  Such correlation learning screens those variables that have
weak marginal correlations with the response.  For classification problems with $Y = \pm 1$, the correlation ranking reduces to selecting features by using two-sample $t$-test statistics.  See Section~\ref{Sec4.2} for additional details.

Other examples of independence learning include methods in microarray data analysis where a two-sample test is used to select significant genes between the treatment and control groups (\cite{DSB03, ST03, FR06, Efron07}), feature ranking using a generalized correlation (\cite{HM09a}),
nonparametric learning under sparse additive models (\cite{RLLW09}), and the method in \cite{HHM08} that uses the marginal bridge estimators for selecting variables in high dimensional sparse regression models.  \cite{HTX09} derive some independence learning rules using tilting methods and empirical likelihood, and propose a bootstrap method to assess the fidelity of feature ranking. In particular, the false discovery rate (FDR) proposed by \cite{BH95} is popularly used in multiple testing for controlling the expected false positive rate. See also \cite{ETST01}, \cite{ABDJ06}, \cite{DJ06}, and \cite{CH09}.

We now discuss the sure screening property of correlation screening.  Let
$\mathcal{M}_* = \{1 \leq j \leq p: \beta_j \neq 0\}$ be the true
underlying sparse model with nonsparsity size $s = |\mathcal{M}_*|$;
the other $p - s$ variables can also be correlated with the
response variable via the link to the predictors
in the true model.   \cite{FL08} consider the case $p \gg n$ with $\log p = O(n^a)$ for some $a \in (0, 1 - 2 \kappa)$, where $\kappa$ is specified below, and Gaussian noise $\varepsilon \sim N(0, \sigma^2)$ for some $\sigma > 0$.  They assume that $\var(Y) = O(1)$,
$\lambda_{\max}(\Sig) = O(n^\tau)$,
\[
\min_{j \in \mathcal{M}_*} |\beta_j| \geq c n^{-\kappa} \quad
\text{and} \quad \min_{j \in \mathcal{M}_*} |\cov(\beta_j^{-1} Y,
X_j)| \geq c,
\]
where $\Sig = \cov(\bx)$, $\kappa, \tau \geq 0$, $c$ is a
positive constant, and the $p$-dimensional covariate vector $\bx$ has an elliptical
distribution with the random matrix $\bX \Sig^{-1/2}$ having a concentration property that holds for Gaussian distributions.  For studies on
the extreme eigenvalues and limiting spectral distributions of large
random matrices, see, e.g., \cite{Silverstein85}, \cite{BY93},
\cite{Bai99}, \cite{Johnstone01}, and \cite{Ledoux01, Ledoux05}.

Under the above regularity conditions, \cite{FL08} show that if
$2 \kappa + \tau < 1$, then there exists some $\theta \in (2 \kappa
+\tau, 1)$ such that when $d_n \sim n^{\theta}$, we have for some
$C > 0$,
\begin{equation} \label{064}
P(\mathcal{M}_* \subset \widehat {\mathcal{M}}_d) = 1 - O(p e^{-C n^{1 - 2
\kappa}/\log n}).
\end{equation}
In particular, this sure screening property entails the sparsity of the model:  $s \leq d_n$. It demonstrates
that SIS can reduce exponentially high dimensionality to
a relatively large scale $d_n \ll  n$, while the reduced model $\widehat{\mathcal{M}}_\gamma$
still contains all the important variables with an overwhelming
probability. In practice, to be conservative we can choose $d = n -
1$ or $[n/\log n]$.  Of course, one can also take final model size $d \geq n$. Clearly larger $d$ means larger probability of
including the true underlying sparse model ${\mathcal{M}}_*$ in the
final model $\widehat{\mathcal{M}}_d$.  See Section 4.3 for further results on sure independence screening.

When the dimensionality is reduced from a large scale $p$ to
a moderate scale $d$ by applying a sure screening method
such as correlation learning, the well-developed variable selection techniques, such as penalized least squares methods, can be applied to the reduced feature space.  This is a powerful tool of SIS based variable selection
methods. The sampling properties of these methods can be easily
obtained by combining the theory of SIS and penalization methods.

\subsection{Feature selection for classification} \label{Sec4.2}
Independence learning has also been widely used for feature selection in high dimensional classification problems. In this section we look at the specific setting of classification and continue the topic of independence learning for variable selection in Section \ref{Sec4.3}. Consider the $p$-dimensional classification between two classes. For
$k \in \{1,2\}$, let $\bX_{k1}$, $\bX_{k2}$, $\cdots$, $\bX_{kn_k}$ be
i.i.d. $p$-dimensional observations from the $k$-th class.
Classification aims at finding a discriminant function $\delta(\bx)$ that classifies new observations as accurately as possible.  The classifier $\delta(\cdot)$ assigns $\bx$ to the class 1 if $\delta(\bx) \geq 0$ and class 2 otherwise.

Many classification methods have been proposed in the literature.  The best classifier is the Fisher $\delta_F(\bx)=(\bx-\bmu)'\bSig^{-1}(\bmu_1-\bmu_2)$
when the data are from the normal distribution with a common covariance matrix: $\bX_{ki} \sim N(\bmu_k, \bSig)$, for $k=1,2$ and $\bmu = (\bmu_1 + \bmu_2)/2$.  However, this method is hard to implement when dimensionality is high due to the difficulty of estimating the unknown covariance matrix $\bSig$.  Hence, the independence rule that involves estimating the diagonal entries of the covariance matrix, with discriminant function
$\delta(\bx)=(\bx-\bmu)'\bD^{-1}(\bmu_1-\bmu_2) $ is frequently
employed for the classification, where $\bD = \diag\{\bSig\}$.  For a survey of recent developments, see \cite{FFW09a}.

Classical methods break down when the dimensionality is
high. As demonstrated by \cite{BL04}, the Fisher discrimination
method no longer performs well in high dimensional settings due to
the diverging spectra and singularity of the sample covariance
matrix. They show that the independence rule overcomes these problems and outperforms the Fisher discriminant in high dimensional setting. However, in practical implementation such as tumor classification using microarray data, one hopes to find tens of genes that have high
discriminative power. The independence rule does not possess the
property of feature selection.

The noise accumulation phenomenon is well-known in the regression setup, but has never been quantified in the classification problem until \cite{FF08}.  They
show that the difficulty of high dimensional
classification is intrinsically caused by the existence of many
noise features that do not contribute to the reduction of
classification error. For example, in linear discriminant analysis
one needs to estimate the class mean vectors and covariance matrix.
Although each parameter can be estimated accurately, aggregated
estimation error can be very large and can significantly
increase the misclassification rate.

Let $\bR_0$ be the common correlation matrix, $\lambda_{\max}(\bR_0)$ be its largest eigenvalue, and $\balpha=\bmu_1-\bmu_2$.  Consider the parameter space
$$
  \Gamma=\{(\balpha,\bSig): \balpha'\bD^{-1}\balpha\geq C_p,\lambda_{\max}(\bR_0)\leq b_0,
\min_{1\leq j\leq p}\sigma_{j}^2>0\} ,
$$
where $C_p$ and $b_0$ are given constants, and $\sigma_j^2$ is the $j$-th diagonal element of $\bSig$. Note that $C_p$ measures the strength of signals.
Let $\hdelta$ be the estimated discriminant function of the
independence rule, obtained by plugging in the sample estimates of
$\balpha$ and $\bD$.  If $\sqrt{n_1n_2/(np)}C_p\rightarrow D_0 \geq 0$, \cite{FF08}
demonstrate that the worst case classification error, $W(\hdelta)$, over the parameter space $\Gamma$ converges:
\begin{equation}\label{e20}
W(\hdelta) \toP 1 - \Phi\Big(\frac{D_0}{2\sqrt{b_0}}\Big),
\end{equation}
where $n = n_1 + n_2$ and $\Phi(\cdot)$ is the cumulative distribution function of the standard normal random variable.

The misclassification rate (\ref{e20}) relates to dimensionality in the term $D_0$, which depends on $C_p/\sqrt{p}$.  This quantifies the tradeoff
between dimensionality $p$ and the overall signal strength $C_p$.
The signal $C_p$ always increases with dimensionality.  If the useful features are located at the first $s$ components, say, then the signals stop increasing when more than $s$ features are used, yet the penalty of using all features is $\sqrt{p}$.  Clearly, using $s$ features can perform much better than using all $p$ features.  The optimal number should be the one that minimizes $C_m/\sqrt{m}$, where the $C_m$ are the signals of the best subset $S$ of $m$ features, defined as $\balpha_S \bD_S^{-1} \balpha_S$, where $\balpha_S$ and $\bD_S$ are the sub-vector and sub-matrix of $\balpha$ and $\bD$ constructed using variables in $S$. The result (\ref{e20}) also indicates that
the independence rule works no better than random
guessing due to noise accumulation, unless the signal levels are
extremely high, say, $\sqrt{n/p}C_p\geq B$ for some $B>0$.  \cite{HPG08b} show that if $C_p^2/p \to \infty$, the classification error goes to zero for a distance-based classifier, which is a specific result of \cite{FF08} with $B = \infty$.

The above results reveal that dimensionality reduction is also very
important for reducing misclassification rate. A popular class of dimensionality reduction techniques is projection. See, for example, principal component analysis in \cite{Ghosh02} and \cite{ZHT04}; partial least squares in
\cite{HP03}, and \cite{Boul04}; and sliced inverse regression in
\cite{CM02}, \cite{ALL03}, and \cite{BP03}.  These projection
methods attempt to find directions that can result in small
classification errors. In fact, the directions that they find usually
put much larger weights on features with large classification power,
which is indeed a type of sparsity in the projection vector.
\cite{FF08} formally show that linear projection methods are likely
to perform poorly unless the projection vector is sparse, namely,
the effective number of selected features is small. This is due to
the aforementioned noise accumulation when estimating $\bmu_1$ and
$\bmu_2$ in high dimensional problems. For
formal results, see Theorem 2 in \cite{FF08}. See also \cite{THNC02}, \cite{DJ08}, \cite{HPS08a}, \cite{HPG08b}, \cite{HC09}, \cite{HM09b}, and \cite{Jin09} for some recent developments in high dimensional classifications.

To select important features,  the two-sample $t$
test is frequently employed (see, e.g., \cite{THNC03}).
The two-sample $t$ statistic for feature $j$ is
\begin{align}\label{e21}
T_j=\frac{\bar{X}_{1j}-\bar{X}_{2j}}{\sqrt{S_{1j}^2/n_1+S_{2j}^2/n_2}},\
j=1,\cdots, p,
\end{align}
where $\bar{X}_{kj}$ and $S_{kj}^2$ are the sample mean and variance
of the $j$-th feature in class  $k$.  This is a specific example of independence learning, which ranks the features according to $|T_j|$.
\cite{FF08} prove that when dimensionality $p$ grows no
faster than the exponential rate of the sample size, if the lowest signal level is
not too small, the two-sample $t$ test can select all important features
with probability tending to 1. Their proof relies on the
deviation results of the two-sample $t$-statistic. See,
e.g., \cite{Hall87, Hall06}, \cite{JSW03}, and \cite{Cao07} for large deviation theory.

Although the $t$ test can correctly select all important features with
probability tending to 1 under some regularity conditions, the
resulting choice is not necessarily optimal, since the noise
accumulation can exceed the signal accumulation for faint features.
Therefore, it is necessary to further single out the most important
features. To address this issue, \cite{FF08} propose the Features
Annealed Independence Rule (FAIR). Instead of constructing
the independence rule using all features, FAIR selects the most
important ones and uses them to construct an independence rule. To
appreciate the idea of FAIR, first note that the relative importance
of features can be measured by $|\alpha_j|/\sigma_j$, where $\alpha_j$ is the $j$-th component of $\balpha = \bmu_1 - \bmu_2$ and $\sigma_j^2$ is the common variance of the $j$-th feature. If such oracle ranking information is available, then one can construct the
independence rule using $m$ features with the largest
$|\alpha_j|/\sigma_j$, with optimal value of $m$ to be determined.
In this case, the oracle classifier takes the form
$$
\hdelta(\bx)=\sum_{j=1}^p\hat{\alpha}_j
(x_j-\hat{\mu}_j)/\hat{\sigma}_j^21_{\{|\alpha_j|/\sigma_j>b\}},
$$
where $b$ is a positive constant. It is easy to see that choosing
the optimal $m$ is equivalent to selecting the optimal $b$. However oracle information is usually unavailable, and one needs
to learn it from the data. Observe that $|\alpha_j|/\sigma_j$ can be
estimated by $|\hat{\alpha}_j|/\hat{\sigma}_j$, where the latter is
in fact $\sqrt{n/(n_1n_2)}|T_j|$, in which the pooled sample variance is used.  This is indeed the same as ranking the feature by using the correlation between the $j$th variable with the class response $\pm 1$ when $n_1 = n_2$ (\cite{FL08}).  Indeed, as pointed out by \cite{HTX08c}, this is always true if the response for the first class is assigned as 1, whereas the response for the second class is assigned as $-n_1/n_2$.  Thus to mimick the oracle, FAIR takes a slightly
different form to adapt to the unknown signal strength
\begin{align}\label{e22}
\fair(\bx)=\sum_{j=1}^p\hat{\alpha}_j (x_j-\hat{\mu}_j)/\hat{\sigma}_j^21_{\{\sqrt{n/(n_1n_2)}|T_j|>b\}}.
\end{align}
It is clear from (\ref{e22}) that FAIR works the same way as if we
first sort the features by the absolute values of their
$t$-statistics in descending order, then take out the first
$m$ features to construct a classifier. The number of features is
selected by minimizing the upper bound of the classification error:
$$
\hat{m}=\arg\max_{1 \leq m\leq
p}\frac{1}{\hat{\lambda}_{\max}^m}\frac{n[\sum_{j=1}^mT_{(j)}^2+m(n_1-n_2)/n]^2}
{mn_1n_2+n_1n_2\sum_{j=1}^mT_{(j)}^2},
$$
where $T_{(1)}^2\geq T_{(2)}^2\geq \cdots \geq T_{(p)}^2$ are the
ordered squared $t$-statistics, and $\hat{\lambda}_{\max}^m$ is the
estimate of the largest eigenvalue of the correlation matrix $\bR_0^m$
of the $m$ most significant features. \cite{FF08} also derive the
misclassification rates of FAIR and demonstrate that it possesses an oracle property.

\subsection{Sure independence screening for generalized linear models} \label{Sec4.3}

Correlation learning cannot be directly applied to the case of discrete covariates such as genetic studies with different genotypes.  The mathematical results and technical arguments in \cite{FL08} rely heavily on the joint normality assumptions.  The natural question is how to screen variables in a more general context, and whether the sure screening property continues to hold with a limited false positive rate.

Consider the generalized linear model (GLIM) with canonical link.  That is, the conditional density is given by
\begin{eqnarray}\label{e23}
f(y |\bx) = \exp \left\{y\theta(\bx) - b(\theta(\bx)) + c(y) \right\},
\end{eqnarray}
for some known functions $b(\cdot)$, $c(\cdot)$, and $\theta(\bx) = \bx\t \bbeta$.  As we consider only variable selection on the mean regression function, we assume without loss of generality that the dispersion parameter $\phi = 1$.  As before, we assume that each variable has been standardized with mean 0 and variance 1.

For GLIM (\ref{e23}), the penalized likelihood (\ref{e02}) is
\begin{equation} \label{e24}
   - n^{-1} \sum_{i=1}^n \ell(\bx_i\t \bbeta, y_i) - \sum_{j=1}^p p_\lambda(|\beta_j|),
\end{equation}
where $\ell(\theta,y) =  b(\theta) - y \theta$.  The maximum marginal likelihood estimator (MMLE) $\hat \bbeta_{j}^M$ is defined as the minimizer of the componentwise regression
\begin{eqnarray} \label{e25}
\hat {\bbeta}_{j}^M = (\hat {\beta}_{j,0}^M, \hat {\beta}_{j}^M )
= \mbox{argmin}_{\beta_0, \beta_j } \sum_{i=1}^n \ell(
\beta_0+\beta_jX_{ij},Y_i),
\end{eqnarray}
where $X_{ij}$ is the $i$th observation of the $j$th variable.
This can be easily computed and its implementation is robust, avoiding numerical instability in ultra-high dimensional problems.  The marginal estimator estimates the wrong object of course, but its magnitude provides useful information for variable screening. \cite{FS09} select a set of variables whose marginal magnitude exceeds a predefined threshold value $\gamma_n$:
\begin{eqnarray}\label{e26}
\widehat {\mathcal{M}}_{\gamma_n} = \{1 \le j \le p: |\hat
\beta_{j}^M| \ge \gamma_n \},
\end{eqnarray}
This is equivalent to ranking features according to the magnitude of MMLEs
$\{|\hat{\beta}_j^M|\}$.
To understand the utility of MMLE, we take the population
version of the minimizer of the componentwise regression to be
\begin{eqnarray*}
{\bbeta}_{j}^M = ({\beta}_{j,0}^M,  {\beta}_{j}^M )\t =
\mbox{argmin}_{\beta_0, \beta_j } {E} \ell( \beta_0+\beta_jX_{j},Y).
\end{eqnarray*}
\cite{FS09} show that $\beta_j^M =0$ if and only if $\cov(X_j, Y) = 0$, and under some additional conditions if $|\cov(X_j, Y)| \ge c_1
n^{-\kappa}$ for $j \in \mathcal{M}_{\star}$, for given positive constants $c_1$ and $\kappa$, then
there exists a  constant $c_2$ such that
\begin{equation}  \label{e27}
  \min_{j \in \mathcal{M}_{\star}}|\beta_j^M| \ge c_2 n^{-\kappa}.
\end{equation}
In words, as long as $X_j$ and $Y$ are somewhat marginally correlated with $\kappa < 1/2$, the marginal signal $\beta_j^M$ is detectable.
They prove further the sure screening property:
\begin{eqnarray} \label{e28}
P\Bigl(\mathcal{M}_{\star} \subset \widehat
{\mathcal{M}}_{\gamma_n}\Bigr)  \to 1
\end{eqnarray}
(the convergence is exponentially fast) if $\gamma_n = c_3 n^{-\kappa}$ with a sufficiently small $c_3$, and that only the size of non-sparse elements (not the dimensionality) matters for the purpose of sure screening property.  For the Gaussian linear model (\ref{e03}) with sub-Gaussian covariate tails, the dimensionality can be as high as $\log p = o(n^{(1-2\kappa)/4})$, a weaker result than that in \cite{FL08} in terms of condition on $p$, but a stronger result in terms of the conditions on the covariates.  For logistic regression with bounded covariates, such as genotypes, the dimensionality can be as high as $\log p = o(n^{1-2\kappa})$.

The sure screening property (\ref{e28}) is only part of the story.  For example, if $\gamma_n = 0$ then all variables are selected and hence (\ref{e28}) holds.  The question is how large the size of the selected model size in (\ref{e26}) with $\gamma_n = c_3 n^{-\kappa}$ should be.  Under some regularity conditions, \cite{FS09} show that with probability tending to one exponentially fast,
\begin{equation} \label{e29}
| \widehat {\cal M}_{\gamma_n}|  = O\{n^{2\kappa}\lambda_{\max} (\bSig)\}.
\end{equation}
In words, the size of selected model depends on how large the thresholding parameter $\gamma_n$ is, and how correlated the features are.
It is of order $O(n^{2 \kappa + \tau})$ if $\lambda_{\max} (\bSig)=O(n^\tau)$.   This is the same or somewhat stronger result than in \cite{FL08} in terms of selected model size, but holds for a much more general class of models.  In particularly, there is no restrictions on $\kappa$ and $\tau$, or more generally $\lambda_{\max} (\bSig)$.

\cite{FS09} also study feature screening by using the marginal likelihood ratio test.  Let $\hat{L}_0 = \mbox{min}_{\beta_0 } n^{-1} \sum_{i=1}^n \ell(\beta_0,Y_i)$ and
\begin{equation} \label{e30}
\hat{L}_j = \hat{L}_0 - \mbox{min}_{\beta_0, \beta_j } n^{-1} \sum_{i=1}^n \ell(
\beta_0+\beta_jX_{ij},Y_i).
\end{equation}
Rank the features according to the marginal utility $\{\hat{L}_j\}$. Thus, select a set of variables
\begin{equation} \label{e31}
\widehat {\mathcal{N}}_{\nu_n} = \{1 \le j \le p_n: \hat{L}_{j} \ge
\nu_n \},
\end{equation}
where $\nu_n$ is a predefined threshold value.  Let $L_j^\star$ be the population counterpart of $\hat{L}_j$.  Then, the minimum signal
$\min_{j \in {\cal M}_*} L_j^\star$ is of order $O(n^{-2 \kappa})$, whereas the individual noise $\hat L_j - L_j^\star = O_p(n^{-1/2})$.  In words, when $\kappa \geq 1/4$, the noise level is larger than the signal.  This is the key technical challenge.  By using the fact that the ranking is invariant to monotonic transformations, \cite{FS09} are able to show that with $\nu_n = c_4 n^{-2 \kappa}$ for a sufficiently small $c_4 > 0$,
$$
P\{{\cal M}_* \subset \widehat {\mathcal{N}}_{\nu_n}, |\widehat {\mathcal{N}}_{\nu_n}| \leq O(n^{2\kappa} \lambda_{\max} (\bSig))\} \to 1.
$$
Thus the sure screening property holds with a limited size of the selected model.

\subsection{Reduction of false positive rate} \label{Sec4.4}

A screening method is usually a crude approach that results in many false positive variables.  A simple idea of reducing the false positive rate is to apply a resampling technique as proposed by \cite{FSW09b}. Split the samples randomly into two halves and let $\hat{\calA}_1$ and $\hat{\cal A}_2$ be the selected sets of active variables based on, respectively, the first half and the second half of the sample.  If $\hat{\cal A}_1$ and $\hat{\cal A}_2$ both have a sure screening property, so does the set $\hat{\cal A}$.  On the other hand, $\hat{\cal A} =\hat{\cal A}_1 \cap \hat{\cal A}_2$ has many fewer falsely selected variables, as an unimportant variable has to be selected twice at random in the ultra-high dimensional space, which is very unlikely.  Therefore, $\hat{\cal A}$ reduces the number of false positive variables.

Write $\mathcal{A}$ for the set of active indices -- that is, the set
containing those indices $j$ for which $\beta_j \neq 0$ in the true
model.   Let $d$ be the size of the selected sets $\calA_1$ and $\calA_2$. Under some exchangeability conditions, \cite{FSW09b} demonstrate that
\begin{equation}  \label{e32}
P(|\widehat{\mathcal{A}} \cap \mathcal{A}^c| \geq r) \leq
\frac{\binom{d}{r}^2}{\binom{p-|\mathcal{A}|}{r}} \leq
\frac{1}{r!}\Bigl(\frac{n^2}{p-|\mathcal{A}|}\Bigr)^r,
\end{equation}
where, for the second inequality, we require that $d \leq n \leq (p - |\mathcal{A}|)^{1/2}$.  In other words, the probability of selecting at least $r$ inactive variables is very small when $n$ is small compared to $p$, such as for the situations discussed in the previous two sections.

\subsection{Iterative sure independence screening} \label{Sec4.5}

SIS uses only the marginal information of the covariates and its
sure screening property can fail when technical conditions
are not satisfied. \cite{FL08} point out three
potential problems with SIS:
\begin{itemize}
\item[a)] {\em (False Negative)} An important predictor that is marginally uncorrelated but jointly correlated with the response cannot be picked by SIS.  An example of this has the covariate vector $\bx$ jointly normal with equi-correlation $\rho$, while $Y$ depends on the covariates through
    $$
       \bx\t \bbeta^\star = X_1 + \cdots + X_J - J \rho X_{J+1}.
    $$
    Clearly, $X_{J+1}$ is independent of $\bx\t \bbeta^\star$ and hence $Y$, yet the regression coefficient $-J \rho$ can be much larger than for other variables.  Such a hidden signature variable cannot be picked by using independence learning, but it has a dominant predictive power on $Y$.

\item[b)] {\em (False Positive)} Unimportant predictors that are highly correlated with the important
predictors can have higher priority to be selected by SIS than
important predictors that are relatively weakly related to the
response.  An illustrative example has
$$
Y = \rho X_0 + X_1 + \cdots + X_J + \varepsilon,
$$
where $X_0$ is independent of the other variables which have a common correlation $\rho$. Then
$\corr(X_j, Y) = J \rho = J \; \corr(X_0, Y)$, for $j = J+1, \cdots, p$,
and $X_0$ has the lowest priority to be selected.

\item[c)] The issue of collinearity among the predictors adds difficulty to the problem of variable selection.
\end{itemize}
Translating a) to microarray data analysis, a two-sample test can never pick up a hidden signature gene.  Yet, missing the hidden signature gene can result in very poor understanding of the molecular mechanism and in poor disease classification.  \cite{FL08} address these issues by proposing an iterative SIS (ISIS) that extends SIS and uses more fully the joint
information of the covariates. ISIS still maintains
computational expediency.

\cite{FSW09b} extend and improve the idea of ISIS from the multiple regression model to the more general loss function (\ref{e15}); this includes, in addition to the log-likelihood, the hinge loss $L(x, y) = (1-xy)_+$ and exponential loss $L(x, y) = \exp(-xy)$ in classification in which $y$ takes values $\pm 1$, among others.  The $\psi$-learning (\cite{STZW03}) can also be cast in this framework. ISIS also allows variable deletion in the process of iteration.  More generally, suppose that our objective is to find a sparse $\bbeta$ to minimize
$$
   n^{-1} \sum_{i=1}^n L(Y_i, \bx_i\t \bbeta) + \sum_{j=1}^p p_\lambda(|\beta_j|).
$$
The algorithm goes as follows.
\begin{enumerate}
\item  Apply an SIS such as (\ref{e31}) to pick a set $\mathcal{A}_1$ of indices of size $k_1$, and then employ a penalized
(pseudo)-likelihood method (\ref{e14}) to select a subset $\mathcal{M}_1$ of these indices.

\item {\em (Large-scale screening)} Instead of computing residuals as in \cite{FL08}, compute
\begin{equation}
L_j^{(2)} = \min_{\beta_0,\sbbeta_{\mathcal{M}_1}, \beta_j} n^{-1}
\sum_{i=1}^n L(Y_i, \beta_0 + \bx_{i, \sM_1}^T \bbeta_{\sM_1} +
X_{ij} \beta_j),  \label{b5}
\end{equation}
for $j \not \in {\mathcal{M}}_1$, where $\bx_{i, \sM_1}$ is the sub-vector of
$\bx_i$ consisting of those elements in ${\mathcal{M}}_1$. This measures the additional contribution of variable $X_j$ in the presence of variables $\bx_{\sM_1}$.  Pick $k_2$ variables with the
smallest $\{L_j^{(2)}, j \not \in {\mathcal{M}}_1\}$ and let $\mathcal{A}_2$ be the resulting set.

\item {\em (Moderate-scale selection)} Use penalized likelihood to obtain
\begin{equation}
\hbbeta_2 = \mbox{argmin}_{\beta_0,\sbbeta_{\sM_1},\sbbeta_{\sA_2}} n^{-1}
\sum_{i=1}^n L(Y_i, \beta_0 + \bx_{i, \sM_1}^T\bbeta_{\sM_1} +
\bx_{i, \sA_2}^T\bbeta_{\sA_2} )+ \sum_{j \in \sM_1 \cup \sA_2}
p_\lambda(|\beta_j|).   \label{e34}
\end{equation}
This gives new active indices
${\mathcal{M}}_2$ consisting of nonvanishing elements of $\hbbeta_2$.
This step also deviates importantly from the approach in \citet{FL08} even in the least squares case. It allows the procedure to delete variables from
the previous selected variables ${\mathcal{M}}_1$.

\item {\em (Iteration)} Iterate the above two steps  until $d$ (a prescribed number) variables are recruited or $\mathcal{M}_\ell = \mathcal{M}_{\ell - 1}$.
\end{enumerate}
The final estimate is then $\hbbeta_{\sM_\ell}$.
In implementation, \cite{FSW09b} choose $k_1 = \lfloor 2d/3 \rfloor$, and
thereafter at the $r$-th iteration, take $k_r = d -
|\mathcal{M}_{r-1}|$.  This ensures that the iterated
versions of SIS take at least two iterations to terminate.
The above method can be considered as an analogue of the least squares ISIS procedure (\cite{FL08}) without explicit definition of the residuals.
\citet{FL08} and \cite{FSW09b} show empirically that the ISIS significantly improves the performance of SIS even in the difficult cases described above.

\section{Sampling properties of penalized least squares} \label{Sec5}
The sampling properties of penalized likelihood estimation (\ref{e02}) have been extensively studied, and a significant amount of work has been contributed to penalized least squares (\ref{e04}). The theoretical studies can be mainly classified into four groups: persistence, consistency and selection consistency, the weak oracle property, and the oracle property (from weak to strong). Again, persistence means consistency of the risk (expected loss) of the estimated model, as opposed to consistency of the estimate of the parameter vector under some loss. Selection consistency means consistency of the selected model. By the weak oracle property, we mean that the estimator enjoys the same sparsity as the oracle estimator with asymptotic probability one, and has consistency. The oracle property is stronger than the weak oracle property in that, in addition to the sparsity in the same sense and consistency, the estimator attains an information bound mimicking that of the oracle estimator. Results have revealed the behavior of different penalty functions and the impact of dimensionality on high dimensional variable selection.

\subsection{Dantzig selector and its asymptotic equivalence to LASSO} \label{Sec5.1}

The $L_1$ regularization (e.g., LASSO) has received much attention due to its convexity and encouraging sparsity solutions. The idea of using the $L_1$ norm can be traced back to the introduction of convex relaxation for deconvolution in \cite{CM73}, \cite{TBM79}, and \cite{SS86}. The use of the $L_1$ penalty has been shown to have close connections to other methods. For example, sparse approximation using an $L_1$ approach is shown in \cite{Girosi98} to be equivalent to support vector machines (\cite{Vapnik95}) for noiseless data. Another example is the asymptotic equivalence between the Dantzig selector (\cite{CT07}) and LASSO.

The $L_1$ regularization has also been used in the Dantzig selector recently proposed by \cite{CT07}, which is defined as the solution to
\begin{equation} \label{e35}
\min \|\bbeta\|_1 \quad \text{subject to} \quad \|n^{-1} \bX\t (\by - \bX \bbeta)\|_\infty \leq \lambda,
\end{equation}
where $\lambda \geq 0$ is a regularization parameter. It was named
after Dantzig because the convex optimization problem (\ref{e35})
can easily be recast as a linear program. Unlike the PLS
(\ref{e04}) which uses the residual sum of squares as a measure of goodness of fit,
the Dantzig selector uses the $L_\infty$ norm of the covariance
vector $n^{-1} \bX\t (\by - \bX \bbeta)$, i.e., the maximum absolute covariance
between a covariate and the residual vector $\by - \bX
\bbeta$, for controlling the model fitting. This $L_\infty$
constraint can be viewed as a relaxation of the normal equation
\begin{equation} \label{e36}
\bX\t \by = \bX\t \bX \bbeta,
\end{equation}
namely, finding the estimator that has the smallest $L_1$-norm in the neighborhood of the least squares estimate.
A prominent feature of the Dantzig selector is its
nonasymptotic oracle inequalities under $L_2$ loss. Consider the
Gaussian linear regression model (\ref{e03}) with $\bveps \sim
N(\bzero, \sigma^2 I_n)$ for some $\sigma > 0$, and assume that each
covariate is standardized to have $L_2$ norm $\sqrt{n}$
(note that we changed the scale of $\bX$ since it was assumed that each
covariate has unit $L_2$ norm in \cite{CT07}). Under the uniform
uncertainty principle (UUP) on the design matrix $\bX$, a condition on the finite condition number for submatrices of $\bX$, they show
that, with high probability, the Dantzig selector $\hbbeta$ mimics
the risk of the oracle estimator up to a logarithmic factor $\log
p$, specifically
\begin{equation} \label{e37}
\|\hbbeta -
\bbeta_0\|_2 \leq C \sqrt{(2 \log p)/n} (\sigma^2 + \sum\nolimits_{j
\in \supp(\bbeta_0)} \beta_{0, j}^2 \wedge \sigma^2)^{1/2},
\end{equation}
where $\bbeta_0 = (\beta_{0, 1}, \cdots, \beta_{0, p})\t$ is the vector of the true regression coefficients, $C$ is some positive constant, and $\lambda \sim \sqrt{(2 \log p)/n}$. Roughly speaking, the UUP condition (see also \cite{DS89, DH01}) requires that all $n \times d$ submatrices of $\bX$ with $d$ comparable to $\|\bbeta_0\|_0$ are uniformly close to orthonormal matrices, which can be stringent in high dimensions. See \cite{FL08} and \cite{CL07} for more discussions. The oracle inequality (\ref{e37}) does not infer much about the sparsity of the estimate.

Shortly after the work on the Dantzig selector, it was observed that the
Dantzig selector and the LASSO share some similarities. \cite{BRT08} present a theoretical comparison of the LASSO and the Dantzig selector in the general high dimensional nonparametric regression model. Under a sparsity scenario, \cite{BRT08} derive parallel oracle inequalities for the prediction risk for both methods, and establish the asymptotic equivalence of the LASSO estimator and the Dantzig selector. More specifically, consider the nonparametric regression model
\begin{equation} \label{e38}
\by = \bff + \bveps,
\end{equation}
where $\bff = (f(\bx_1), \cdots, f(\bx_n))\t$ with $f$ an unknown $p$-variate function, and $\by$, $\bX = (\bx_1, \cdots, \bx_n)\t$, and $\bveps$ are the same as in (\ref{e03}). Let $\{f_1,\cdots, f_M\}$ be a finite dictionary of $p$-variate functions. As pointed out in \cite{BRT08}, $f_j$'s can be a collection of basis functions for approximating $f$, or estimators arising from $M$ different methods. For any $\bbeta = (\beta_1,\cdots, \beta_M)\t$, define $f_{\bbeta}=\sum_{j=1}^M\beta_jf_j$. Then similarly to (\ref{e04}) and (\ref{e35}), the LASSO estimator $\widehat{f}_L$ and Dantzig selector $\widehat{f}_D$ can be defined accordingly as $f_{\hbbeta_L}$ and $f_{\hbbeta_D}$ with $\hbbeta_L$ and $\hbbeta_D$ the corresponding $M$-vectors of minimizers. In both formations, the empirical norm $\|f_j\|_n = \sqrt{n^{-1}\sum_{i=1}^n f_j^2(\bx_i)}$ of $f_j$ is incorporated as its scale. \cite{BRT08} show that under the restricted eigenvalue condition on the Gram matrix and some other regularity conditions, with significant probability, the difference between
$\|\widehat{f}_D - f\|_n^2$ and $\|\widehat{f}_L - f\|_n^2$
is bounded by a product of three factors. The first factor $s \sigma^2/n$ corresponds to the prediction error rate in regression with $s$ parameters, and the other two factors including $\log M$ reflect the impact of a large number of regressors. They further prove sparsity oracle inequalities for the prediction loss of both estimators. These inequalities entail that the distance between the prediction losses of the Dantzig selector and the LASSO estimator is of the same order as the distances between them and their oracle approximations.

\cite{BRT08} also consider the specific case of a linear model (\ref{e03}), say (\ref{e38}) with true regression function $\bff = \bX \bbeta_0$. If $\bveps \sim
N(\bzero, \sigma^2 I_n)$ and some regularity conditions hold, they show that, with large probability, the $L_q$ estimation loss for $1\leq q \leq 2$ of the Dantzig selector $\hbbeta_D$ is simultaneously given by
\begin{equation} \label{e39}
\|\hbbeta_D - \bbeta_0\|_q^q \leq C \sigma^q \left(1+\sqrt{s/m}\right)^{2(q-1)}s\left(\frac{\log p}{n}\right)^{q/2} ,
\end{equation}
where $s= \|\bbeta_0\|_0$, $m \geq s$ is associated with the strong restricted eigenvalue condition on the design matrix $\bX$, and $C$ is some positive constant. When $q = 1$, they prove (\ref{e39}) under a (weak) restricted eigenvalue condition that does not involve $m$. \cite{BRT08} also derive similar inequalities to (\ref{e39}) with slightly different constants on the $L_q$ estimation loss, for $1\leq q \leq 2$, of the LASSO estimator $\hbbeta_L$. These results demonstrate the approximate equivalence of the Dantzig selector and the LASSO. The similarity between the Dantzig selector and LASSO has also been discussed in \cite{EHT07}. \cite{Lounici08} derives the $L_\infty$ convergence rate and studies a sign concentration property simultaneously for the LASSO estimator and the Dantzig selector under a mutual coherence condition.

Note that the covariance vector $n^{-1} \bX\t (\by - \bX \bbeta)$ in the formulation of Dantzig selector (\ref{e35}) is exactly the negative gradient of $(2n)^{-1} \|\by - \bX \bbeta\|^2$ in PLS
(\ref{e04}). This in fact entails that the Dantzig selector and  the LASSO estimator are identical under some suitable conditions, provided that the same regularization parameter $\lambda$ is used in both methods. For example, \cite{MRY07} give a diagonal dominance condition of the $p \times p$ matrix $(\bX\t \bX)^{-1}$ that ensures
their equivalence. This condition implicitly assumes $p \leq n$. \cite{JRL09}
present a formal necessary and sufficient condition, as well as easily verifiable sufficient conditions ensuring the identical solution of the Dantzig selector and the LASSO estimator when the
dimensionality $p$ can exceed sample size $n$.

\subsection{Model selection consistency of LASSO} \label{Sec5.2}

There is a huge literature devoted to studying the statistical properties of LASSO and related methods. This $L_1$ method as well as its variants have also been extensively studied in such other areas as compressed sensing. For example, \cite{GR04} show that under some regularity conditions the LASSO-type procedures are persistent under quadratic loss for dimensionality of polynomial growth, and \cite{Greenshtein06} extends the results to more general loss functions. \cite{Meinshausen07} presents similar results for the LASSO for dimensionality of exponential growth and finite nonsparsity size, but its persistency rate is slower than that of a relaxed LASSO. For consistency and selection consistency results see \cite{DET06}, \cite{MB06}, \cite{Wainwright06}, \cite{ZY06}, \cite{BTW07}, \cite{BRT08}, \cite{VDG08}, and \cite{ZH08}, among others.

As mentioned in the previous section, consistency results for the LASSO hold under some conditions on the design matrix. For the purpose of variable selection, we are also concerned with the sparsity of the estimator, particularly its model selection consistency meaning that the estimator $\hbbeta$ has the same support as the true regression coefficients vector $\bbeta_0$ with asymptotic probability one. \cite{ZY06} characterize the model selection consistency of the LASSO by studying a stronger but technically more convenient property of sign consistency: $P(\sgn(\hbbeta) = \sgn(\bbeta_0)) \rightarrow 1$ as $n \rightarrow \infty$. They show that the weak irrepresentable condition
\begin{equation} \label{e40}
\|\bX_2\t \bX_1 (\bX_1\t \bX_1)^{-1} \sgn(\bbeta_1)\|_\infty < 1
\end{equation}
is necessary for sign consistency of the LASSO, and the strong irrepresentable condition, which requires that the left-hand side of (\ref{e40}) be uniformly bounded by a positive constant $C < 1$, is sufficient for sign consistency of the LASSO, where $\bbeta_1$ is the subvector of $\bbeta_0$ on its support $\supp(\bbeta_0)$, and $\bX_1$ and $\bX_2$ denote the submatrices of the $n \times p$ design matrix $\bX$ formed by columns in $\supp(\bbeta_0)$ and its complement, respectively. See also \cite{Zou06} for the fixed $p$ case. However, the irrepresentable condition can become restrictive in high dimensions. See Section \ref{Sec5.4} for a simple illustrative example, because the same condition shows up in a related problem of sparse recovery by using $L_1$ regularization. This demonstrates that in high dimensions, the LASSO estimator can easily select an inconsistent model, which explains why the LASSO tends to include many false positive variables in the selected model.

To establish the weak oracle property of the LASSO, in addition to the sparsity characterized above, we need its consistency. To this end, we usually need the condition on the design matrix that
\begin{equation} \label{e41}
\|\bX_2\t \bX_1 (\bX_1\t \bX_1)^{-1}\|_\infty \leq C
\end{equation}
for some positive constant $C < 1$, which is stronger than the strong irrepresentable condition. It says that the $L_1$-norm of the regression coefficients of each inactive variable regressed on $s$ active variables must be uniformly bounded by $C < 1$.  This shows that the capacity of the LASSO for selecting a consistent model is very limited, noticing also that the $L_1$-norm of the regression coefficients typically increase with $s$.  See, e.g., \cite{Wainwright06}. As discussed above, condition (\ref{e41}) is a stringent condition in high dimensions for the LASSO estimator to enjoy the weak oracle property. The model selection consistency of the LASSO in the context of graphical models has been studied by \cite{MB06}, who consider Gaussian graphical models with polynomially growing numbers of nodes.

\subsection{Oracle property} \label{Sec5.3}
What are the sampling properties of penalized least squares (\ref{e04}) and penalized likelihood estimation (\ref{e02}) when the penalty function $p_\lambda$ is no longer convex? The oracle property (\cite{FL01}) provides a nice conceptual framework for understanding the statistical properties of high dimensional variable selection methods.

In a seminal paper, \cite{FL01} build the theoretical foundation of nonconvex penalized least squares or, more generally, nonconcave penalized likelihood for variable selection. They introduce the oracle property for model selection. An estimator $\hbbeta$ is said to have the oracle property if it enjoys sparsity in the sense that $\hbbeta_2 = \bzero$ with probability tending to 1 as $n \rightarrow \infty$, and $\hbbeta_1$ attains an
information bound mimicking that of the oracle estimator, where $\hbbeta_1$ and $\hbbeta_2$ are  the subvectors of $\hbbeta$ formed by components in $\supp(\bbeta_0)$ and $\supp(\bbeta_0)^c$, respectively, while the oracle knows the true model $\supp(\bbeta_0)$ beforehand. The oracle properties of penalized least squares estimators can be understood in the more general framework of penalized likelihood estimation. \cite{FL01} study the oracle properties of nonconcave penalized likelihood estimators in the finite-dimensional setting, and \cite{FP04} extend their results to the moderate dimensional setting with $p = o(n^{1/5})$ or $o(n^{1/3})$.

More specifically, without loss of generality, assume that the true regression coefficients vector is $\bbeta_0 = (\bbeta_1\t, \bbeta_2\t)\t$ with $\bbeta_1$ and $\bbeta_2$ the subvectors of nonsparse and sparse elements respectively:  $\|\bbeta_1\|_0 = \|\bbeta_0\|_0 $ and $\bbeta_2 = \bzero$. Let
$a_n = \|p_\lambda'(|\bbeta_1|)\|_\infty$ and $b_n = \|p_\lambda''(|\bbeta_1|)\|_\infty$.
\cite{FL01} and \cite{FP04} show that, as long as $a_n, b_n = o(1)$, under some regularity conditions there exists a local maximizer $\hbbeta$ to the penalized likelihood (\ref{e02}) such that
\begin{equation} \label{e042}
\|\hbbeta - \bbeta_0\|_2 = O_P(\sqrt{p} (n^{-1/2} + a_n)).
\end{equation}
This entails that choosing the regularization parameter $\lambda$ with $a_n = O(n^{-1/2})$ gives a root-$(n/p)$ consistent
penalized likelihood estimator. In particular, this is the case when the SCAD penalty is used if $\lambda = o(\min_{1 \leq j \leq s}|\beta_{0, j}|)$, where $\bbeta_1 = (\beta_{0, 1}, \cdots, \beta_{0, s})\t$.  Recently, \cite{FL09} gave a sufficient condition under which the solution is unique.

\cite{FL01} and \cite{FP04} further prove the oracle properties
of penalized likelihood estimators under some additional
regularity conditions. Let
$\bSig = \diag\{p_\lambda''(|\bbeta_1|)\}$
and
$
    \bar{p}_{\lambda}(\bbeta_1) = \sgn(\bbeta_1) \circ p_\lambda'(|\bbeta_1|),
$
where $\circ$ denotes the the Hadamard (componentwise) product.
Assume that $\lambda = o(\min_{1 \leq j \leq s}|\beta_{0, j}|)$, $\sqrt{n/p} \lambda \rightarrow \infty$ as $n \rightarrow \infty$, and the penalty function $p_{\lambda}$ satisfies
$\liminf_{n \rightarrow \infty} \liminf_{t \rightarrow 0+} p'_{\lambda}(t)/\lambda > 0$.
They show that if $p = o(n^{1/5})$, then with probability tending to 1 as $n \rightarrow \infty$, the root-$(n/p)$ consistent local maximizer $\hbbeta = (\hbbeta_1\t, \hbbeta_2\t)\t$ satisfies the following
\begin{itemize}
\item[a)] (Sparsity) $\hbbeta_2 = \bzero$;

\item[b)] (Asymptotic normality)
\begin{equation} \label{e043}
\sqrt{n} \bA_n \bI_1^{-1/2} (\bI_1 + \bSig) [\hbbeta_1 - \bbeta_1 + (\bI_1 + \bSig)^{-1} \bar{p}_{\lambda}(\bbeta_1)] \toD N(\bzero, \bG),
\end{equation}
\end{itemize}
where $\bA_n$ is a $q \times s$ matrix such that $\bA_n \bA_n\t \rightarrow \bG$, a $q \times q$ symmetric positive definite matrix, $\bI_1 = \bI(\bbeta_1)$ is the Fisher information matrix knowing the true model $\supp(\bbeta_0)$, and $\hbbeta_1$ is a subvector of $\hbbeta$ formed by components in $\supp(\bbeta_0)$.

Consider a few penalties. For the SCAD penalty, the condition $\lambda = o(\min|\bbeta_1|)$ entails that both $\bar{p}_{\lambda}(\bbeta_1)$ and $\Sigma$ vanish asymptotically. Therefore, the asymptotic normality (\ref{e043}) becomes
\begin{equation} \label{e044}
\sqrt{n} \bA_n \bI_1^{1/2} (\hbbeta_1 - \bbeta_1) \toD N(\bzero, \bG),
\end{equation}
which shows that $\hbbeta_1$ has the same asymptotic efficiency as the MLE of $\bbeta_1$ knowing the true model in advance. This demonstrates that the resulting penalized likelihood estimator is as efficient as the oracle one. For the $L_1$ penalty (LASSO), the root-$(n/p)$ consistency of $\hbbeta$ requires $\lambda = a_n = O(n^{-1/2})$, whereas the oracle property requires $\sqrt{n/p} \lambda \rightarrow \infty$ as $n \rightarrow \infty$. However, these two conditions are incompatible, which suggests that the LASSO estimator
generally does not have the oracle property. This is intrinsically due to the fact that the $L_1$ penalty does not satisfy the unbiasedness condition.

It has indeed been shown in \cite{Zou06} that the LASSO estimator does not have the oracle property even in the finite parameter setting. To address the bias issue of LASSO, he proposes the adaptive LASSO by using an adaptively weighted $L_1$ penalty. More specifically, the weight vector is $|\hbbeta|^{-\gamma}$ for some $\gamma > 0$ with the power understood componentwise, where $\hbbeta$ is an initial root-$n$ consistent estimator of $\bbeta_0$. Since $\hbbeta$ is root-$n$ consistent, the constructed weights can separate important variables from unimportant ones.   This is an attempt to introduce the SCAD-like penalty to reduce the biases.  From (\ref{e12}), it can easily be seen that the adaptive LASSO is just a specific solution to penalized least squares using LLA.  As a consequence, \cite{Zou06} shows that the adaptive LASSO has the oracle property under some regularity conditions. See also \cite{ZH08}.

\subsection{Additional properties of SCAD estimator}

In addition to the oracle properties outlined in the last section and also in Section 6.2, \cite{KCO08} and \cite{KK09} provide insights into the SCAD estimator. They attempt to answer the question of when the oracle estimator $\hat{\bbeta}^o$ is a local minimizer of the penalized least squares with the SCAD penalty, when the SCAD estimator and the oracle estimator coincide, and how to check whether a local minimizer is a global minimizer.  The first two results are indeed stronger than the oracle property as they show that the SCAD estimator is the oracle estimator itself rather than just mimicking its performance.

Recall that all covariates have been standardized.  The follow assumption is needed.
\begin{itemize}
\item []{\bf Condition A}:  The nonsparsity size is $s_n = O(n^{c_1})$ for some $0 < c_1 < 1$, the minimum eignvalue of the correlation matrix of those active variables is bounded away from zero, and the minimum signal $\min_{1 \leq j \leq s_n} |\beta_{j}| > c_3 n^{-(1-c_2)/2}$ for some constant $c_2 \in (c_1, 1]$.
\end{itemize}
Under Condition A, \cite{KCO08} prove that if $E \varepsilon_i^{2k} < \infty$ for the linear model (\ref{e03}),
\begin{equation} \label{eq45a}
       P(\mbox{$\hat{\bbeta}^o$ is a local minima of PLS with the SCAD penalty}) \to 1,
\end{equation}
provided that $\lambda_n = o(n^{-\{1-(c_2-c_1)\}/2})$ and $p_n = o \{ (\sqrt{n} \lambda_n)^{2k}\}$.  This shows that the SCAD method produces the oracle estimator.  When $k$ is sufficiently large, the dimensionality $p_n$ can be of any polynomial order. For the Gaussian error, the result holds even with NP-dimensionality.  More precisely,  for the Gaussian errors, they show that (\ref{eq45a}) holds for $p_n = O(\exp(c_4 n))$ and $\lambda_n = O(n^{-(1-c_5)/2})$, where $0 < c_4 < c_5 < c_2 - c_1$.  The question then arises naturally whether the global minimizer of penalized least squares with the SCAD penalty is the oracle estimator.  \cite{KCO08} give an affirmative answer:
with probability tending to one, the global minimizer of penalized least squares with the SCAD penalty is the same as the oracle estimator when the correlation matrix of all covariates is bounded away from zero and infinity (necessarily, $p_n \leq n$).

\cite{KK09} also give a simple condition under which the SCAD estimator is unique and is a global minimizer (see also the simple conditions in \cite{FL09} for a more general problem).  They also provide sufficient conditions to check whether a local minimizer is a global minimizer.  They show that the SCAD method produces the oracle estimator,
$$
  P\{\mbox{The SCAD estimator} = \hat{\bbeta}^o\} \to 1,
$$
under conditions similar to Condition A, even when the minimum eigenvalue of the correlation matrix of all variables converges to zero.

\subsection{Sparse recovery and compressed sensing} \label{Sec5.4}

Penalized $L_1$ methods have been widely applied in areas including compressed sensing (\cite{Donoho06a}). In those applications, we want to find good sparse representations or approximations of signals that can greatly improve efficiency of data storage and transmission. We have no intent here to survey results on compressed sensing. Rather, we would like to make some innate connections of the problem of sparse recovery in the noiseless case to model selection. Unveiling the role of penalty functions in sparse recovery can give a simplified view of the role of penalty functions in high dimensional variable selection as the noise level approaches zero. In particular, we see that concave penalties are advantageous in sparse recovery, which is in line with the advocation of folded concave penalties for variable selection as in \cite{FL01}.

Consider the noiseless case $\by = \bX \bbeta_0$ of the linear model (\ref{e03}). The problem of sparse recovery aims to find the sparsest possible solution
\begin{equation} \label{e45}
\arg \min \|\bbeta\|_0 \quad \text{subject to} \quad \by = \bX \bbeta.
\end{equation}
The solution to $\by = \bX \bbeta$ is not unique when the $n \times p$ matrix $\bX$ has rank less than $p$, e.g., when $p > n$. See \cite{DE03} for a characterization of the identifiability of the minimum $L_0$ solution $\bbeta_0$. Although by its nature, the $L_0$ penalty is the target penalty for sparse recovery, its computational complexity makes it infeasible to implement in high dimensions. This motivated the use of penalties that are computationally tractable
relaxations or approximations to the $L_0$ penalty. In particular, the convex $L_1$ penalty provides a nice convex relaxation and has attracted much attention. For properties of various $L_1$ and related methods see, for example, the Basis Pursuit in \cite{CDS99}, \cite{DE03},
\cite{Donoho04}, \cite{Fuchs04}, \cite{CT05, CT06},
\cite{DET06}, \cite{Tropp06},
\cite{CWB08}, and \cite{CXZ09}.

More generally, we can replace the $L_0$ penalty in (\ref{e45}) by a penalty function $\rho(\cdot)$ and consider the $\rho$-regularization problem
\begin{equation} \label{e46}
\min \sum_{j = 1}^p \rho(|\beta_j|) \quad \text{subject to} \quad \by = \bX \bbeta.
\end{equation}
This constrained optimization problem is closely related to the PLS in
(\ref{e04}). A great deal of research has contributed to identifying conditions on $\bX$ and $\bbeta_0$ that ensure the $L_1/L_0$ equivalence, i.e., the $L_1$-regularization (\ref{e46}) gives the same solution $\bbeta_0$. For example, \cite{Donoho04} contains deep results and shows that the individual equivalence of $L_1$/$L_0$ depends only on $\supp(\bbeta_0)$ and $\bbeta_{0}$ on its support. See also \cite{DH01} and \cite{Donoho06b}. In a recent work, \cite{LF09} present a sufficient condition that ensures the $\rho/L_0$ equivalence for concave penalties. They consider increasing and concave penalty functions $\rho(\cdot)$ with finite maximum concavity (curvature). The convex $L_1$ penalty falls at the boundary of this class of penalty functions. Under these regularity conditions, they show that $\bbeta_0$ is a local minimizer of (\ref{e46}) if there exists some $\epsilon \in (0, \min_{j \leq s} |\beta_{0, j}|)$ such that
\begin{equation} \label{e47}
\max_{\bu \in \mathcal{U}_\epsilon} \|\bX_2\t \bX_1 (\bX_1\t \bX_1)^{-1} \bu\|_\infty < \rho'(0+),
\end{equation}
where $\mathcal{U}_\epsilon = \{\sgn(\bbeta_{1}) \circ  \rho'(|\bv|): \|\bv - \bbeta_1\|_\infty \leq \epsilon\}$, the notation being that of the previous two sections.

When the $L_1$ penalty is used,
$\mathcal{U}_\epsilon$ contains a single point $\sgn(\bbeta_1)$ with $\sgn$ understood componentwise.
In this case, condition (\ref{e47}) becomes the weak irrepresentable condition (\ref{e40}). In fact the $L_1/L_0$ equivalence holds provided that (\ref{e40}) weakened to nonstrict inequality is satisfied. However, this condition can become restrictive in high dimensions. To appreciate this, look at an example given in \cite{LF09}. Suppose that $\bX_1 = (\bx_1, \cdots, \bx_s)$ is orthonormal, $\by = \sum_{j = 1}^s \beta_{0, j} \bx_j$ with $|\beta_{0, 1}| = \cdots = |\beta_{0, s}|$, $\bx_{s + 1}$ has unit $L_2$ norm and correlation $r$ with $\by$, and all the rest of the $\bx_j$'s are orthogonal to $\{\bx_j\}_{j=1}^s$. The above condition becomes $|r| \leq s^{-1/2}$. This demonstrates that the $L_1$ penalty can fail to recover the sparsest solution $\bbeta_0$ when the maximum correlation of the noise variable and response is moderately high, which, as explained in the Introduction, can easily happen in high dimensions.

On the other hand, the concavity of the penalty function $\rho$ entails that its derivative $\rho'(t)$ is deceasing in $t \in [0, \infty)$. Therefore, condition (\ref{e47}) can be (much) less restrictive for concave penalties other than $L_1$. This shows the advantage of concave penalties in sparse recovery, which is consistent with similar understandings in variable selection in \cite{FL01}.

\section{Oracle property of penalized likelihood with ultra-high dimensionality} \label{Sec6}

As shown in Section \ref{Sec4}, large-scale screening and moderate-scale selection is a good strategy for variable selection in ultra-high dimensional feature spaces.  A less stringent screening (i.e., a larger selected model size in (\ref{e23})) will have a higher probability of retaining all important variables. It is important to study the limits of the dimensionality that nonconcave penalized likelihood methods can handle. The existing result of \cite{FP04} is too weak in terms of the dimensionality allowed for high dimensional modeling;  they deal with too broad a class of models.

What are the roles of the dimensionality $p$ and nonsparsity size $s$?
What is the role of penalty functions? Does the oracle property continue to hold in ultra-high dimensional feature spaces? These questions have been driving the theoretical development of high dimensional variable selection. For example, \cite{Koltchinskii08} obtains oracle inequalities for penalized least squares with entropy penalization, and \cite{VDG08} establishes a
nonasymptotic oracle inequality for the Lasso estimator as the empirical risk minimizer in high dimensional generalized linear models. There are relatively few studies on the statistical properties of high dimensional variable selection methods by using regularization with nonconvex penalties. More recent studies on this topic include \cite{HHM08}, \cite{KCO08}, \cite{MvB08}, \cite{LF09}, \cite{Zhang09}, and \cite{FL09}, among others.

\subsection{Weak oracle property} \label{Sec6.1}

An important step towards the understanding of the oracle property of penalized likelihood methods in ultra-high dimensions is the weak oracle property for model selection, introduced by \cite{LF09} in the context of penalized least squares. An estimator $\hbbeta$ is said to have the weak oracle property if it is uniformly consistent and enjoys sparsity in the sense of $\hbbeta_2 = \bzero$ with probability tending to 1, i.e. model selection consistency, where $\hbbeta_2$ is the subvector of $\hbbeta$ formed by components in $\supp(\bbeta_0)^c$ and the oracle knows the true model $\supp(\bbeta_0)$ beforehand. This property is weaker than the oracle property in \cite{FL01}. Consistency is derived under $L_\infty$ loss, mainly due to the technical difficulty of proving the existence of a solution to the nonlinear equation that characterizes the nonconcave penalized likelihood estimator. It is important to study the rate of the probability bound for sparsity and the rate of convergence for consistency. The dimensionality $p$ usually enters the former rate explicitly, from which we can see the allowable growth rate of $p$ with sample size $n$.

Consider the PLS problem (\ref{e04}) with penalty function $p_\lambda$.  Let $\rho(t; \lambda) = \lambda^{-1} p_\lambda(t)$ and write it as $\rho(t)$ whenever there is no confusion.  \cite{LF09} and \cite{FL09} consider the following class of penalty functions:
\begin{itemize}
\item $\rho(t; \lambda)$ is increasing and concave in $t \in [0, \infty)$, and has a continuous derivative $\rho'(t; \lambda)$ with $\rho'(0+; \lambda) > 0$. In addition, $\rho'(t; \lambda)$ is increasing in $\lambda \in (0, \infty)$ and $\rho'(0+; \lambda)$ is independent of $\lambda$.
\end{itemize}
This is a wide class of concave penalties including SCAD and MCP, and the $L_1$ penalty at its boundary. \cite{LF09} establish a nonasymptotic weak oracle property for the PLS estimator. They consider (\ref{e03}) with $\bveps \sim N(\bzero, \sigma^2 I_n)$. The notation here is the same as in Section \ref{Sec5}. Assume that each column of the $n \times p$ design matrix $\bX$ (covariate) is standardized to have $L_2$ norm $\sqrt{n}$ (or of this order), and let
$d_n = 2^{-1} \min\left\{\left|\beta_{0, j}\right|: \beta_{0, j} \neq 0\right\}$
be the minimal signal. The following condition is imposed on the design matrix
\begin{equation} \label{e48}
\|\bX_2\t \bX_1 (\bX_1\t \bX_1)^{-1}\|_\infty \leq \min (  C \frac{\rho'(0+)}{\rho'(d_n)}, O(n^{\alpha_1}))
\end{equation}
where $\alpha_1 \geq 0$, $C \in (0, 1)$, and $\rho$ is associated with the regularization parameter $\lambda \sim n^{\alpha_1 - 1/2} u_n$. Here $\{u_n\}$ is a sequence of positive numbers diverging to infinity. Clearly, for the $L_1$ penalty, condition (\ref{e48}) becomes (\ref{e41}) which is a somewhat stronger form of the strong irrepresentable condition in \cite{ZY06}. Condition (\ref{e48}) consists of two parts:  the first part is intrinsic to the penalty function whereas the second part is purely a technical condition. For folded-concave penalties other than $L_1$, the intrinsic condition is much more relaxed: the intrinsic upper bound is $C < 1$ for the $L_1$ penalty whereas it is $\infty$ when $d_n \gg \lambda_n$ for the SCAD type of penalty.  In other words, the capacity for LASSO to have model selection consistency is limited, independent of model signals, whereas no limit is imposed for SCAD type of penalties when the signals are strong enough.  In general, the concavity of $\rho(\cdot)$ guarantees condition (\ref{e48}) and is more relaxed than the $L_1$ penalty.

Under the above and some additional regularity conditions, if $ \|(n^{-1} \bX_1\t \bX_1)^{-1}\|_\infty = O(n^{\alpha_0})$ for some $\alpha_0 \geq 0$,
\cite{LF09} show that for sufficiently large $n$, with probability at least $1 - \frac{2}{\sqrt{\pi}} p u_n^{-1} e^{-u_n^2/2}$, the PLS estimator $\hbbeta = (\hbbeta_1\t, \hbbeta_2\t)\t$ satisfies the following
\begin{itemize}
\item[a)] (Sparsity) $\hbbeta_2 = \bzero$;

\item[b)] ($L_\infty$ loss) $\|\hbbeta_1 - \bbeta_1\|_\infty = O(n^{\alpha_0 - 1/2} u_n)$,
\end{itemize}
where $\hbbeta_1$ is a subvector of $\hbbeta$ formed by components in $\supp(\bbeta_0)$.  In particular, when the signals are so sparse that $s$ is finite, $\alpha_0 = 0$ for all non-degenerate problems.   In this case, by taking $u_n^2 = c \log p$ for $c \geq 2$ so that the probability $1 - \frac{2}{\sqrt{\pi}} p u_n^{-1} e^{-u_n^2/2} \to 1$, we have
$\|\hbbeta_1 - \bbeta_1\|_\infty = O_P(n^{-1/2} \sqrt{\log p})$.
As an easy consequence of the general result,
\begin{equation} \label{e49}
\|\hbbeta - \bbeta_0\|_2 = O_P(\sqrt{s} n^{\alpha_0 - 1/2} u_n)
\end{equation}
when $p = o(u_n e^{u_n^2/2})$. The dimensionality $p$ is allowed to grow up to exponentially fast with $u_n$. More specifically, $u_n$ can be allowed as large as $o(n^{1/2 - \alpha_0 - \alpha_1} d_n)$ and thus $\log p = o(n^{1 - 2(\alpha_0 + \alpha_1)} d_n^2)$.
This shows that a weaker minimal signal needs slower growth of dimensionality for successful variable selection. From their studies, we also see the known fact that concave penalties can reduce the biases of estimates.

Recently, \cite{FL09} extended the results of \cite{LF09} and established a nonasymptotic weak oracle property for non-concave penalized likelihood estimator in generalized linear models with ultra-high dimensionality. In their weak oracle property, they relax the term $u_n$ from the consistency rate. A similar condition to (\ref{e48}) appears, which shows the drawback of the $L_1$ penalty. The dimensionality $p$ is allowed to grow at a non-polynomial (NP) rate. Therefore, penalized likelihood methods can still enjoy the weak oracle property in ultra-high dimensional space.

\subsection{Oracle property with NP-dimensionality} \label{Sec6.2}

A long-standing question is whether the penalized likelihood methods enjoy the oracle property (\cite{FL01}) in ultra-high dimensions.  This issue has recently been addressed by \cite{FL09} in the context of generalized linear models. Such models include the commonly used linear, logistic, and Poisson regression models.

More specifically \cite{FL09} show that, under some regularity conditions, there exists a local maximizer $\hbbeta = (\hbbeta_1\t, \hbbeta_2\t)\t$ of the penalized likelihood (\ref{e02}) such that $\hbbeta_2 = \bzero$ with probability tending to 1 and $\|\hbbeta - \bbeta_0\|_2 = O_P(\sqrt{s} n^{-1/2})$, where $\hbbeta_1$ is a subvector of $\hbbeta$ formed by components in $\supp(\bbeta_0)$ and $s = \|\bbeta_0\|_0$. They further establish asymptotic normality and thus the oracle property. The conditions are less restrictive for such concave penalties as SCAD. In particular, their results suggest that the $L_1$ penalized likelihood estimator generally cannot achieve the consistent rate of $O_P(\sqrt{s} n^{-1/2})$ and does not have the oracle property when the dimensionality $p$ is diverging with the sample size $n$. This is consistent with results in \cite{FL01}, \cite{FP04}, and \cite{Zou06}.

It is natural to ask when the non-concave penalized likelihood estimator is also a global maximizer of the penalized likelihood (\ref{e02}). \cite{FL09} give characterizations of such a property from two perspectives: global optimality and restricted global optimality. In particular, they show that under some regularity conditions, the SCAD penalized likelihood estimator can be identical to the oracle estimator. This feature of the SCAD penalty is not shared by the $L_1$ penalty.

\section{Concluding remarks} \label{Sec7}

We now have a better picture of the role of penalty functions and the impact of dimensionality on high dimensional regression and classification. The whole story of high dimensional statistical learning is far from complete. New innovative techniques are needed and critical analyses of their relative merits are required. Issues include the characterization of optimality properties, the selection of data-driven penalty functions and parameters, the confidence in selected models and estimated parameters, group variable selection and its properties, inference after model selection, the incorporation of information on covariates, nonparametric statistical learning, manifold learning, compressed sensing, developments of high dimensional statistical techniques in other important statistical contexts, and development of robust and user-friendly algorithms and software.  High dimensional statistical learning is developed to confront and address the challenges in the frontiers of scientific research and technological innovation. It interfaces nicely with many scientific disciplines and will undoubtedly further advances on emerging societal needs.

\begin{singlespace}

\end{singlespace}

\begin{thebibliography}{999}
\bibitem[Abramovich et al.(2006)]{ABDJ06}
Abramovich, F., Benjamini, Y., Donoho, D., and Johnstone, I. (2006). Adapting to unknown
sparsity by controlling the false discovery rate. \textit{Ann. Statist.} \textbf{34}, 584--653.

\bibitem[Akaike(1973)]{Akaike73}
Akaike, H. (1973). Information theory and an extension of the maximum likelihood principle. In \textit{Second International Symposium on Information Theory} (eds. B. N. Petrov and F. Csaki), 267--281, Akademiai Kiado, Budapest.

\bibitem[Akaike(1974)]{Akaike74}
Akaike, H. (1974). A new look at the statistical model identification. \textit{IEEE Trans. Auto. Control} \textbf{19}, 716--723.

\bibitem [Antoniadis(1996)]{Antoniadis96}
Antoniadis, A. (1996). Smoothing noisy data with tapered coiflets series.   {\em Scand. J. Statist.} {\bf  23}, 313--330.

\bibitem[Antoniadis and Fan(2001)]{AF01}
Antoniadis, A. and Fan, J. (2001). Regularization of wavelets approximations (with discussion). \textit{J. Amer. Statist. Assoc.} \textbf{96}, 939--967.

\bibitem[Antoniadis, Lambert-Lacroix and Leblanc(2003)]{ALL03}
Antoniadis, A., Lambert-Lacroix, S., and Leblanc, F. (2003).
Effective dimension reduction methods for tumor classification using
gene expression data. \textit{Bioinformatics} \textbf{19}, 563--570.

\bibitem[Bai(1999)]{Bai99}
Bai, Z. D. (1999). Methodologies in spectral analysis of large dimensional random matrices, a review. \textit{Statist. Sin.} \textbf{9}, 611--677.

\bibitem[Bai and Yin(1993)]{BY93}
Bai, Z. D. and Yin, Y. Q. (1993). Limit of smallest eigenvalue of a large dimensional sample covariance matrix. \textit{Ann. Probab.} \textbf{21}, 1275--1294.

\bibitem[Barron, Birge and Massart(1999)]{BBM99}
Barron, A., Birge, L., and Massart, P. (1999). Risk bounds for model selection via penalization. {\em Probab. Theory Related Fields} {\bf 113}, 301--413.

\bibitem[Benjamini and Hochberg(1995)]{BH95}
Benjamini, Y. and Hochberg, Y. (1995). Controlling the false discovery rate: a practical and
powerful approach to multiple testing. \textit{J. Roy. Statist.
Soc. B} \textbf{57}, 289--300.

\bibitem[Bickel(1983)]{Bickel83}
Bickel, P. J. (1983). Minimax estimation of a normal mean subject to doing well at a point.
In \textit{Recent Advances in Statistics} (ed. by M. H. Rizvi, J. S. Rustagi, and D. Siegmund), 511--528, Academic Press, NewYork.

\bibitem[Bickel(2008)]{Bickel08}
Bickel, P. J. (2008).  Discussion of ``Sure independence screening
for ultrahigh dimensional feature space''. \textit{J. Roy. Statist.
Soc. B} \textbf{70}, 883--884.

\bibitem[Bickel and Levina(2004)]{BL04}
Bickel, P. J. and Levina, E. (2004). Some theory for
Fisher's linear discriminant function, ``naive Bayes", and some
alternatives when there are many more variables than observations.
\textit{Bernoulli} \textbf{10}, 989--1010.

\bibitem[Bickel and Li(2006)]{BL06}
Bickel, P. J. and Li, B. (2006). Regularization in statistics (with discussion). \textit{Test} \textbf{15}, 271--344.

\bibitem[Bickel, Ritov and Tsybakov(2008)]{BRT08}
Bickel, P. J., Ritov, Y., and Tsybakov, A. (2008). Simultaneous analysis of LASSO and Dantzig selector. \textit{Ann. Statist.}, to appear.

\bibitem[Boulesteix(2004)]{Boul04}
Boulesteix, A. (2004). PLS Dimension reduction for classification with microarray
data. \textit{Statistical Applications in Genetics and Molecular
Biology} \textbf{3}, 1--33.

\bibitem[Bradic, Fan and Wang(2009)]{BFW09}
Bradic, J., Fan, J., and Wang, W. (2009). Penalized composite quasi-likelihood for high-dimensional variable selection. {\em Manuscript}.

\bibitem[Breiman(1995)]{Breiman95}
Breiman, L. (1995). Better subset regression using the non-negative garrote. \textit{Technometrics} \textbf{37}, 373--384.

\bibitem[Breiman(1996)]{Breiman96}
Breiman, L. (1996). Heuristics of instability and stabilization in model selection. \textit{Ann. Statist.} \textbf{24}, 2350--2383.

\bibitem[Bunea, Tsybakov and Wegkamp(2007)]{BTW07}
Bunea, F., Tsybakov, A., and Wegkamp, M. H. (2007). Sparsity oracle inequalities for the LASSO. \textit{Elec. Jour. Statist.} \textbf{1}, 169--194.

\bibitem[Bura and Pfeiffer(2003)]{BP03}
Bura, E. and Pfeiffer, R. M. (2003).
Graphical methods for class prediction using
dimension reduction techniques on DNA microarray data. \textit{Bioinformatics} \textbf{19}, 1252--1258.

\bibitem[Cai and Lv(2007)]{CL07}
Cai, T. and Lv, J. (2007). Discussion: The Dantzig selector: statistical estimation when $p$ is much larger than $n$. \textit{Ann. Statist.} \textbf{35}, 2365--2369.

\bibitem[Cai, Xu and Zhang(2009)]{CXZ09}
Cai, T., Xu, G., and Zhang, J. (2009). On recovery of sparse signals via $\ell_1$ minimization.
\textit{IEEE Transactions on Information Theory}, to appear.

\bibitem[Candes and Tao(2005)]{CT05}
Candes, E. J. and Tao, T. (2005). Decoding by linear programming. \textit{IEEE Trans. Inform. Theory} \textbf{51}, 4203--4215.

\bibitem[Candes and Tao(2006)]{CT06}
Candes, E. J. and Tao, T. (2006). Near-optimal signal recovery from random projections: Universal encoding strategies? \textit{IEEE Trans. Inform. Theory} \textbf{52}, 5406--5425.

\bibitem[Candes and Tao(2007)]{CT07}
Candes, E. and Tao, T. (2007). The Dantzig selector: Statistical estimation when $p$ is much larger than $n$ (with discussion). \textit{Ann. Statist.} \textbf{35}, 2313--2404.

\bibitem[Cand\`{e}s, Wakin and Boyd(2008)]{CWB08}
Cand\`{e}s, E. J., Wakin, M. B., and Boyd, S. P. (2008). Enhancing sparsity by reweighted $\ell_1$ minimization. \textit{J. Fourier Anal. Appl.} \textbf{14}, 877--905.

\bibitem[Cao (2007)]{Cao07}
Cao, H.Y. (2007). Moderate deviations for two sample $t$-statistics. {\em ESAIM Probab. Statist.} {\bf 11}, 264--627.

\bibitem[Chen, Donoho and Saunders(1999)]{CDS99}
Chen, S., Donoho, D. L., and Saunders, M. (1999). Atomic decomposition by basis pursuit. \textit{SIAM Journal of Scientific Computing} \textbf{20}, 33--61.

\bibitem[Chiaromonte and Martinelli(2002)]{CM02}
Chiaromonte, F. and Martinelli, J. (2002). Dimension reduction strategies for
analyzing global gene expression data with a response. \textit{Mathematical Biosciences} \textbf{176}, 123--144.

\bibitem[Claerbout and Muir(1973)]{CM73}
Claerbout, J. F. and Muir, F. (1973). Robust modeling of erratic data. \textit{Geophysics}
\textbf{38}, 826--844.

\bibitem[Clarke and Hall(2009)]{CH09}
Clarke, S. and Hall, P. (2009). Robustness of multiple testing procedures against dependence. \textit{Ann. Statist.} \textbf{37}, 332--358.

\bibitem[Daubechies, Defrise and De Mol(2004)]{DDD04}
Daubechies, I., Defrise, M., and De Mol, C. (2004).
    An iterative thresholding algorithm for linear inverse
    problems with a sparsity constraint. {\em Comm. Pure Appl. Math.}
    {\bf 57},  1413--1457.

\bibitem[Dempster, Laird and Rubin(1977)]{Dempster77}
Dempster, A. P., Laird, N. M., and Rubin, D. B. (1977). Maximum likelihood from incomplete
data via the EM algorithm. \textit{J. Royal Statist. Soc. B} \textbf{39}, 1--38.

\bibitem[Donoho(2000)]{Donoho00}
Donoho, D. L. (2000). High-dimensional data analysis: The curses and blessings of dimensionality. \textit{Aide-Memoire of a Lecture at AMS Conference on Math Challenges of the 21st Century}.

\bibitem[Donoho(2004)]{Donoho04}
Donoho, D. L. (2004). Neighborly polytopes and sparse solution of underdetermined linear equations. \textit{Technical Report}, Department of Statistics, Stanford University.

\bibitem[Donoho(2006a)]{Donoho06a}
Donoho, D. L. (2006a). Compressed sensing. \textit{IEEE Trans. Inform. Theory} \textbf{52}, 1289--1306.

\bibitem[Donoho(2006b)]{Donoho06b}
Donoho, D. L. (2006b). For most large undetermined systems of linear equations the minimal $\ell_1$-norm solution is the sparsest solution. \textit{Comm. Pure Appl. Math.} \textbf{59}, 797--829.

\bibitem[Donoho and Elad(2003)]{DE03}
Donoho, D. L. and Elad, M. (2003). Optimally sparse representation in general (nonorthogonal)
dictionaries via $\ell_1$ minimization. \textit{Proc. Natl. Acad. Sci.} \textbf{100}, 2197--2202.

\bibitem[Donoho, Elad and Temlyakov(2006)]{DET06}
Donoho, D. L., Elad, M., and Temlyakov, V. (2006). Stable recovery of sparse overcomplete representations in the presence of noise. \textit{IEEE Trans. Inform. Theory} \textbf{52}, 6--18.

\bibitem[Donoho and Huo(2001)]{DH01}
Donoho, D. and Huo, X. (2001). Uncertainty principles and ideal atomic decomposition. \textit{IEEE Trans. Inform. Theory} \textbf{47}, 2845--2862.

\bibitem[Donoho and Jin(2006)]{DJ06}
Donoho, D. and Jin, J. (2006). Asymptotic minimaxity of false discovery rate thresholding
for sparse exponential data. \textit{Ann. Statist.} \textbf{34}, 2980--3018.

\bibitem[Donoho and Jin(2008)]{DJ08}
Donoho, D. and Jin, J. (2008). Higher criticism thresholding: optimal feature selection when
useful features are rare and weak. \textit{Proc. Natl. Acad. Sci.} \textbf{105}, 14790--14795.

\bibitem[Donoho and Johnstone(1994)]{DJ94}
Donoho, D. L. and Johnstone, I. M. (1994). Ideal spatial adaptation by wavelet shrinkage.
\textit{Biometrika} \textbf{81}, 425--455.

\bibitem[Donoho and Stark(1989)]{DS89}
Donoho, D. L. and Stark, P. B. (1989). Uncertainty principles and signal recovery.
\textit{SIAM Journal on Applied Mathematics} \textbf{49}, 906--931.

\bibitem[Dudoit, Shaffer and Boldrick(2003)]{DSB03}
Dudoit, S., Shaffer, J. P., and Boldrick, J. C. (2003).
    Multiple hypothesis testing in microarray experiments.
    {\em Statist. Sci.} {\bf 18}, 71--103.

\bibitem[Efron(2007)]{Efron07}
Efron, B. (2007). Correlation and large-scale simultaneous
        significance testing.   {\em
        Jour. Amer. Statist. Assoc.} {\bf 102}, 93--103.

\bibitem[Efron et al.(2004)]{EHJT04}
Efron, B., Hastie, T., Johnstone, I., and Tibshirani, R. (2004). Least angle regression (with discussion). \textit{Ann. Statist.} \textbf{32}, 407--499.

\bibitem[Efron, Hastie and Tibshirani(2007)]{EHT07}
Efron, B., Hastie, T., and Tibshirani, R. (2007). Discussion: The Dantzig selector: statistical estimation when $p$ is much larger than $n$. \textit{Ann. Statist.} \textbf{35}, 2358--2364.

\bibitem[Efron et al.(2001)]{ETST01}
Efron, B., Tibshirani, R., Storey, J., and Tusher, V. (2001). Empirical Bayes analysis of a
microarray analysis experiment. \textit{J. Amer. Statist. Assoc.} \textbf{99}, 96--104.

\bibitem[Fan(1997)]{Fan97}
Fan, J. (1997). Comments on ``Wavelets in statistics: A review" by A. Antoniadis. \textit{J. Italian Statist. Assoc.} \textbf{6}, 131--138.

\bibitem[Fan and Fan(2008)]{FF08}
Fan, J. and Fan, Y. (2008). High-dimensional classification using features annealed independence
rules. \textit{Ann. Statist.} \textbf{36}, 2605--2637.

\bibitem[Fan, Fan and Lv(2008)]{FFL08}
Fan, J., Fan, Y., and Lv, J. (2008). High dimensional covariance matrix estimation using a factor model. \textit{Journal of Econometrics} \textbf{147}, 186--197.

\bibitem[Fan, Fan and Wu(2010)]{FFW09a}
Fan, J., Fan, Y., and Wu, Y. (2010). {\em High dimensional classification}.
To appear in {\em High-dimensional Statistical Inference} (T. Cai and X. Shen, eds.), World Scientific, New Jersey.

\bibitem[Fan and Li(2001)]{FL01}
Fan, J. and Li, R. (2001). Variable selection via nonconcave
penalized likelihood and its oracle properties. \textit{J. Amer.
Statist. Assoc.} \textbf{96}, 1348--1360.

\bibitem[Fan and Li(2006)]{FL06}
Fan, J. and Li, R. (2006). Statistical challenges with high dimensionality: Feature selection in knowledge discovery. \textit{Proceedings of the International Congress of Mathematicians} (M. Sanz-Sole, J. Soria, J.L. Varona, and J. Verdera, eds.), Vol. \textbf{III}, 595--622.

\bibitem[Fan and Lv(2008)]{FL08}
Fan, J. and Lv, J. (2008). Sure independence screening for ultrahigh dimensional feature space (with discussion). \textit{J. Roy. Statist. Soc. B} \textbf{70}, 849--911.

\bibitem[Fan and Lv(2009)]{FL09}
Fan, J. and Lv, J. (2009). Properties of non-concave penalized likelihood with NP-dimensionality. \textit{Manuscript}.

\bibitem[Fan and Peng(2004)]{FP04}
Fan, J. and Peng, H. (2004). Nonconcave penalized likelihood
with diverging number of parameters. \textit{Ann. Statist.}
\textbf{32}, 928--961.

\bibitem[Fan and Ren(2006)]{FR06}
Fan, J. and Ren, Y. (2006). Statistical analysis of DNA microarray data. \textit{Clinical
Cancer Research} \textbf{12}, 4469--4473.

\bibitem[Fan, Samworth and Wu(2009)]{FSW09b}
Fan, J., Samworth, R., and Wu, Y. (2009). Ultrahigh dimensional variable selection: beyond the linear model.  \textit{Journal of Machine Learning Research} {\bf 10}, 1829--1853.

\bibitem[Fan and Song(2009)]{FS09}
Fan, J. and Song, R. (2009). Sure independence screening in generalized
linear models with NP-dimensionality. Revised for \textit{Ann. Statist.}

\bibitem[Foster and George(1994)]{FG94}
Foster, D. and George, E. (1994). The risk inflation criterion for multiple regression. \textit{Ann. Statist.} \textbf{22}, 1947--1975.

\bibitem[Frank and Friedman(1993)]{FF93}
Frank, I. E. and Friedman, J. H. (1993). A statistical view of some chemometrics regression tools (with discussion). \textit{Technometrics} \textbf{35}, 109--148.

\bibitem[Friedman et al.(2007)]{FHHT07}
Friedman, J., Hastie, T., H\"{o}fling, H., and Tibshirani, R. (2007). Pathwise coordinate optimization. \textit{Ann. Appl. Statist.} \textbf{1}, 302--332.

\bibitem[Friedman, Hastie and Tibshirani(2007)]{FHT07}
Friedman, J., Hastie, T., and Tibshirani, R. (2007). Sparse inverse covariance estimation with the LASSO. \textit{Manuscript}.

\bibitem[Fu(1998)]{Fu98}
Fu, W. J. (1998). Penalized regression: the bridge versus the LASSO. {\em Journal of Computational and Graphical Statistics} {\bf 7}, 397--416.

\bibitem[Fuchs(2004)]{Fuchs04}
Fuchs, J.-J. (2004). Recovery of exact sparse representations in the presence of noise. \textit{Proc. IEEE Int. Conf. Acoustics, Speech, and Signal Processing}, 533--536.

\bibitem[Ghosh(2002)]{Ghosh02}
Ghosh, D. (2002). Singular value decomposition regression
modeling for classification of tumors from microarray experiments.
\textit{Proceedings of the Pacific Symposium on Biocomputing},
11462--11467.

\bibitem[Girosi(1998)]{Girosi98}
Girosi, F. (1998). An equivalence between sparse approximation and support
vector machines. \textit{Neural Comput.} \textbf{10}, 1455--1480.

\bibitem[Greenshtein(2006)]{Greenshtein06}
Greenshtein, E. (2006). Best subset selection, persistence in high-dimensional statistical learning and optimization under $\ell_1$ constraint. \textit{Ann. Statist.} \textbf{34}, 2367--2386.

\bibitem[Greenshtein and Ritov(2004)]{GR04}
Greenshtein, E. and Ritov, Y. (2004). Persistence in high-dimensional predictor
selection and the virtue of overparametrization. \textit{Bernoulli} \textbf{10}, 971--988.

\bibitem[Hall(1987)]{Hall87}
Hall, P. (1987). Edgeworth expansion for Student's $t$ statistic under minimal moment conditions.
\textit{Ann. Probab.} \textbf{15}, 920--931.

\bibitem[Hall(2006)]{Hall06}
Hall, P. (2006). Some contemporary problems in statistical sciences. \textit{The Madrid Intelligencer}, to appear.

\bibitem[Hall and Chan(2009)]{HC09}
Hall, P. and Chan, Y.-B. (2009). Scale adjustments for classifiers in high-dimensional, low sample size settings. \textit{Biometrika} \textbf{96}, 469--478.

\bibitem[Hall, Marron and Neeman(2005)]{HMN05}
Hall, P., Marron, J. S., and Neeman, A. (2005). Geometric representation of high dimension, low sample size data. \textit{J. R. Statist. Soc. B} \textbf{67}, 427--444.

\bibitem[Hall and Miller(2009a)]{HM09a}
Hall, P. and Miller, H. (2009a).  Using generalized correlation to
        effect variable selection in very high dimensional problems.  {\em Jour. Comput. Graphical. Statist.}, to appear.

\bibitem[Hall and Miller(2009b)]{HM09b}
Hall, P. and Miller, H. (2009b). Recursive methods for variable selection in very high dimensional classification. \textit{Manuscript}.

\bibitem[Hall, Park and Samworth(2008)]{HPS08a}
Hall, P., Park, B., and Samworth, R. (2008). Choice of neighbor order in nearest-neighbor classification. \textit{Ann. Statist.} \textbf{5}, 2135-2152.

\bibitem[Hall, Pittelkow and Ghosh(2008)]{HPG08b}
Hall, P., Pittelkow, Y., and Ghosh, M. (2008). Theoretic measures of relative performance of
classifiers for high dimensional data with small sample sizes. \textit{J. Roy. Statist.
    Soc. B} \textbf{70}, 158--173.

\bibitem[Hall, Titterington and Xue(2008)]{HTX08c}
    Hall, P., Titterington, D. M., and Xue, J.-H. (2008).
    Discussion of ``Sure independence screening
    for ultrahigh dimensional feature space''. \textit{J. Roy. Statist.
    Soc. Ser. B} \textbf{70}, 889--890.

\bibitem[Hall, Titterington and Xue(2009)]{HTX09}
    Hall, P., Titterington, D. M., and Xue, J.-H. (2009).  Tilting methods for assessing the influence of components in a classifier.  {\em
    Jour. Roy. Statist. Soc.  B}, to appear.

\bibitem[Hastie, Tibshirani and Friedman(2009)]{HTF09}
Hastie, T., Tibshirani, R., and Friedman, J. (2009). \textit{The Elements of Statistical Learning: Data Mining, Inference, and Prediction} (2nd edition). Springer-Verlag, New York.

\bibitem[Huang, Horowitz and Ma(2008)]{HHM08}
Huang, J., Horowitz, J., and Ma, S. (2008). Asymptotic properties of bridge estimators in sparse high-dimensional regression models. \textit{Ann. Statist.} \textbf{36}, 587--613.

\bibitem[Huang and Pan(2003)]{HP03}
Huang, X. and Pan, W. (2003). Linear regression and two-class classification with
gene expression data. \textit{Bioinformatics} \textbf{19}, 2072--2078.

\bibitem[Hunter and Lange(2000)]{HL00}
Hunter, D. R. and Lange, K. (2000). Rejoinder to discussion of ``Optimization transfer
using surrogate objective functions.'' \textit{J. Comput. Graph. Statist.} \textbf{9}, 52--59.

\bibitem[Hunter and Li(2005)]{HL05}
Hunter, D. R. and Li, R. (2005). Variable selection using MM algorithms. \textit{Ann. Statist.} \textbf{33}, 1617--1642.

\bibitem[James, Radchenko and Lv(2009)]{JRL09}
James, G., Radchenko, P., and Lv, J. (2009). DASSO: connections between the Dantzig selector and LASSO. \textit{J. Roy. Statist. Soc. B} \textbf{71}, 127--142.

\bibitem[Jin(2009)]{Jin09}
Jin, J. (2009). Impossibility of successful classification when useful features are rare and
weak. \textit{Proc. Natl. Acad. Sci.} \textbf{106}, 8859--8864.

\bibitem[Jing, Shao and Wang(2003)]{JSW03}
Jing, B. Y., Shao, Q.-M., and Wang, Q. Y. (2003). Self-normalized Cram\'{e}r type large deviations
for independent random variables. \textit{Ann. Probab.} \textbf{31}, 2167--2215.

\bibitem[Johnstone(2001)]{Johnstone01}
Johnstone, I. M. (2001). On the distribution of the largest eigenvalue in principal components analysis. \textit{Ann. Statist.} \textbf{29}, 295--327.

\bibitem[Kim, Choi and Oh(2008)]{KCO08}
Kim, Y., Choi, H., and Oh, H.S. (2008). Smoothly clipped absolute deviation on high dimensions.  {\em Jour. Amer. Statist. Assoc.} {\bf 103}, 1665--1673.

\bibitem[Kim and Kwon(2009)]{KK09} Kim, Y. and Kwon, S. (2009).  On the global optimum of the SCAD penalized estimator.  {\em Manuscript}.


\bibitem[Koenker(1984)]{Koenker84}
Koenker, R. (1984).  A note on $L$-estimates for linear models.
      {\em Statistics and Probability Letters} {\bf 2}, 323--325.

\bibitem [Koenker and Bassett(1978)]{KB78}
       Koenker, R. and Bassett, G. (1978).  Regression quantiles.
      {\em Econometrica} {\bf 46}, 33--50.

\bibitem[Koltchinskii(2008)]{Koltchinskii08}
Koltchinskii, V. (2008). Sparse recovery in convex hulls via entropy penalization. \textit{Ann. Statist.}, to appear.

\bibitem[Lange(1995)]{Lange95}
Lange, K. (1995). A gradient algorithm locally equivalent to the EM algorithm. \textit{J. Royal Statist. Soc. B} \textbf{57}, 425--437.

\bibitem[Ledoux(2001)]{Ledoux01}
Ledoux, M. (2001). \textit{The Concentration of Measure Phenomenon}. American Mathematical Society, Cambridge.

\bibitem[Ledoux(2005)]{Ledoux05}
Ledoux, M. (2005). Deviation inequalities on largest eigenvalues. \textit{Manuscript}.

\bibitem[Lounici(2008)]{Lounici08}
Lounici, K. (2008). Sup-norm convergence rate and sign concentration property of Lasso and Dantzig estimators. \textit{Electronic Journal of Statistics} \textbf{2}, 90--102.

\bibitem[Lv and Fan(2009)]{LF09}
Lv, J. and Fan, Y. (2009). A unified approach to model selection and sparse recovery using regularized least squares. \textit{Ann. Statist.} \textbf{37}, 3498--3528.

\bibitem[Lv and Liu(2008)]{LL08}
Lv, J. and Liu, J. S. (2008). New principles for model selection when models are possibly
misspecified. \textit{Manuscript}.

\bibitem[Mallows(1973)]{Mallows73}
Mallows, C. L. (1973). Some comments on $C_p$. \textit{Technometrics} \textbf{15}, 661--675.

\bibitem[McCullagh and Nelder(1989)]{MN89}
McCullagh, P. and Nelder, J. A. (1989). \textit{Generalized Linear Models}. Chapman and Hall, London.

\bibitem[Meier, van de Geer and B\"{u}hlmann(2008)]{MvB08}
Meier, L., van de Geer, S., and B\"{u}hlmann, P. (2008). The group LASSO for logistic regression. \textit{J. R. Statist. Soc. B} \textbf{70}, 53--71.

\bibitem[Meinshausen(2007)]{Meinshausen07}
Meinshausen, N. (2007). Relaxed LASSO. \textit{Computnl Statist. Data Anal.} \textbf{52}, 374--393.

\bibitem[Meinshausen and B\"{u}hlmann(2006)]{MB06}
Meinshausen, N. and B\"{u}hlmann, P. (2006). High dimensional graphs and variable selection with the LASSO. \textit{Ann. Statist.} \textbf{34}, 1436--1462.

\bibitem[Meinshausen, Rocha and Yu(2007)]{MRY07}
Meinshausen, N., Rocha, G., and Yu, B. (2007). Discussion: A tale of three cousins: LASSO, L2Boosting and Dantzig. \textit{Ann. Statist.} \textbf{35}, 2373--2384.

\bibitem[Osborne, Presnell and Turlach(2000)]{OPT00}
Osborne, M. R., Presnell, B., and Turlach, B. A. (2000). On the LASSO and its dual. {\em Journal of Computational and Graphical Statistics} {\bf  9}, 319--337.

\bibitem[Ravikumar et al.(2009)]{RLLW09}
Ravikumar, P., Lafferty, J., Liu, H., and Wasserman, L. (2009). Sparse additive models. \textit{Jour. Roy. Statist. Soc. B}, to appear.

\bibitem[Rosset and Zhu(2007)]{RZ07}
Rosset, S. and Zhu, J. (2007).  Piecewise linear regularized solution paths. {\em Ann. Statist.} {\bf 35}, 1012--1030.

\bibitem[Santosa and Symes(1986)]{SS86}
Santosa, F. and Symes, W. W. (1986). Linear inversion of band-limited reflection
seismograms. \textit{SIAM J. Sci. Statist. Comput.} \textbf{7}, 1307--1330.

\bibitem[Schwartz(1978)]{Schwartz78}
Schwartz, G. (1978). Estimating the dimension of a model. \textit{Ann. Statist.} \textbf{6}, 461--464.

\bibitem [Shen et al.(2003)]{STZW03}
Shen, X., Tseng, G.C., Zhang, X., and Wong, W.H. (2003). On $\psi$-learning. {\em Jour. Ameri. Statist. Assoc.} {\bf 98}, 724--734.

\bibitem[Silverstein(1985)]{Silverstein85}
Silverstein, J. W. (1985). The smallest eigenvalue of a large dimensional Wishart matrix. \textit{Ann. Probab.} \textbf{13}, 1364--1368.

\bibitem[Storey and Tibshirani(2003)]{ST03}
Storey, J. D. and Tibshirani R. (2003).  Statistical
    significance for genome-wide studies. {\em Proc. Natl. Aca. Sci.}
     {\bf 100}, 9440--9445.

\bibitem[Taylor, Banks and McCoy(1979)]{TBM79}
Taylor, H. L., Banks, S. C., and McCoy, J. F. (1979). Deconvolution with the $\ell_1$ norm. \textit{Geophysics} \textbf{44}, 39--52.

\bibitem[Tibshirani(1996)]{Tibshirani96}
Tibshirani, R. (1996). Regression shrinkage and selection via the LASSO. \textit{J. Roy. Statist. Soc. B} \textbf{58}, 267--288.

\bibitem[Tibshirani et al.(2002)]{THNC02}
Tibshirani, R., Hastie, T., Narasimhan, B., and Chu, G. (2002). Diagnosis of multiple cancer
types by shrunken centroid of gene expression. \textit{Proc. Natl. Acad. Sci.} \textbf{99}, 6567--6572.

\bibitem[Tibshirani et al.(2003)]{THNC03}
Tibshirani, R., Hastie, T., Narasimhan, B., and Chu, G. (2003). Class prediction by
nearest shrunken centroids, with applications to DNA microarrays. \textit{Statist.
Sci.} \textbf{18}, 104--117.

\bibitem[Tropp(2006)]{Tropp06}
Tropp, J. A. (2006). Just relax: convex programming methods for identifying sparse signals in noise. \textit{IEEE Transactions on Information Theory} \textbf{5}, 1030--1051.

\bibitem[van de Geer(2008)]{VDG08}
van de Geer, S. (2008). High-dimensional generalized linear models and the LASSO. \textit{Ann. Statist.} \textbf{36}, 614--645.

\bibitem[Vapnik(1995)]{Vapnik95}
Vapnik, V. (1995). \textit{The Nature of Statistical Learning}. Springer-Verlag, New York.

\bibitem[Wainwright(2006)]{Wainwright06}
Wainwright, M. J. (2006). Sharp thresholds for high-dimensional and noisy recovery of sparsity. \textit{Technical Report}, Department of Statistics, UC Berkeley.

\bibitem[Wang, Li and Tsai(2007)]{WLT07}
Wang, H., Li, R., and Tsai, C.-L. (2007). Tuning parameter selectors for the
smoothly clipped absolute deviation method. {\em Biometrika} {\bf 94}, 553--568.

\bibitem[Wu and Lang(2008)]{WL08}
Wu, T. T. and Lange, K. (2008). Coordinate descent algorithms
     for LASSO penalized regression.  {\em Ann. Appl. Stat.}
     {\bf 2}, 224--244.

\bibitem[Yuan and Lin(2006)]{YL06}
Yuan, M. and Lin, Y. (2006). Model selection and estimation in regression with grouped variables. {\em Jour. Roy. Statist.
    Soc. B} {\bf 68}, 49--67.

\bibitem[Zhang(2009)]{Zhang09}
Zhang, C.-H. (2009). Penalized linear unbiased selection. \textit{Ann. Statist.}, to appear.

\bibitem[Zhang and Huang(2008)]{ZH08}
Zhang, C.-H. and Huang, J. (2008). The sparsity and bias of the LASSO selection in high-dimensional linear regression. \textit{Ann. Statist.} \textbf{36}, 1567--1594.

\bibitem[Zhang and Li(2009)]{ZL09}
Zhang, Y. and Li, R. (2009). Iterative conditional maximization algorithm for nonconcave
penalized likelihood. \textit{IMS Lecture Notes-Monograph Series}, to appear.

\bibitem[Zhao and Yu(2006)]{ZY06}
Zhao, P. and Yu, B. (2006). On model selection consistency of LASSO. \textit{Journal of Machine Learning Research} \textbf{7}, 2541--2563.

\bibitem[Zou(2006)]{Zou06}
Zou, H. (2006). The adaptive LASSO and its oracle properties. \textit{J. Amer. Statist. Assoc.} \textbf{101}, 1418--1429.

\bibitem[Zou and Hastie(2005)]{ZH05}
Zou, H. and  Hastie, T. (2005). Regularization and variable selection
    via the elastic net. {\em Jour. Roy. Statist. Soc. B} {\bf 67}, 301--320.

\bibitem[Zou, Hastie and Tibshirani(2004)]{ZHT04}
Zou, H., Hastie, T., and Tibshirani. R. (2004). Sparse principal component analysis.
\textit{Technical report}.

\bibitem[Zou and Li(2008)]{ZL08}
Zou, H. and Li, R. (2008). One-step sparse estimates in
nonconcave penalized likelihood models (with discussion).
\textit{Ann. Statist.} \textbf{36}, 1509--1566.

\bibitem[Zou and Yuan(2008)]{ZY08}
Zou, H. and Yuan, M. (2008).  Composite quantile regression
    and the oracle model selection theory. {\em Ann. Statist.} {\bf 36}, 1108--1126.
\end{thebibliography}
\end{document}